\documentclass{article}
\usepackage[utf8]{inputenc}
\usepackage{amsmath, amssymb, amsthm}
\usepackage{xcolor}
\usepackage{subcaption}
\usepackage[font=footnotesize]{caption}
\usepackage{graphics, graphicx} 
\usepackage{cleveref}
\newcommand{\M}{\mathcal{M}}
\newcommand{\so}{\mathfrak{so}}
\newcommand{\se}{\mathfrak{se}}

\newcommand{\dexp}{\mathrm{dexp}}

\newcommand{\g}{\mathfrak{g}}
\newcommand{\Ad}{\mathrm{Ad}}
\newcommand{\ad}{\mathrm{ad}}

\newcommand{\tol}{\mathrm{tol}}
\newcommand{\infgen}{\ensuremath{\psi_*}}
\newcommand{\aux}[1]{\tilde{#1}}

\definecolor{BurntOrange}{rgb}{0.8, 0.33, 0.0}
\newtheorem{example}{Example}

\newcommand{\RE}[1]{#1}

\title{Lie Group integrators for mechanical systems}
\author{Elena Celledoni \and Ergys Çokaj \and Andrea Leone \and Davide Murari \and Brynjulf Owren}
\date{Revised version July 2021, last compiled: \today}

\begin{document}

\maketitle

\begin{abstract}
 Since they were introduced in the 1990s, Lie group integrators have become a  method of choice in many application areas. These include multibody dynamics, shape analysis, data science, image registration and biophysical simulations. Two important classes of intrinsic Lie group integrators are the Runge--Kutta--Munthe--Kaas methods and the commutator free Lie group integrators.
 We give a short introduction to these classes of methods. The Hamiltonian framework is attractive for many mechanical problems, and in particular we shall consider Lie group integrators for problems on cotangent bundles of Lie groups where a number of different formulations are possible. There is a natural symplectic structure on such manifolds and through variational principles one may derive symplectic Lie group integrators. We also consider the practical aspects of the implementation of Lie group integrators, such as adaptive time stepping. The theory is illustrated by applying the methods to two nontrivial applications in mechanics. One is the N-fold spherical pendulum where we introduce the restriction of the adjoint action of the group $SE(3)$ to $TS^2$, the tangent bundle of the two-dimensional sphere. Finally, we show how Lie group integrators can be applied to model the controlled path of a payload being transported by two rotors. This problem is modeled on $\mathbb{R}^6\times \left(SO(3)\times \mathfrak{so}(3)\right)^2\times (TS^2)^2$ and put in a format where Lie group integrators can be applied.
\end{abstract}
\section{Introduction}

In many physical problems, including multi-body dynamics, the configuration space is not a linear space, but rather consists of a collection of rotations and translations.
A simple example is the free rigid body whose configuration space consists of rotations in 3D. A more advanced example is the simplified model of the human body, where the skeleton at a given time is described as a system of interacting rods and joints.
Mathematically, the structure of such problems is usually best described as a manifold. Since manifolds by definition can be equipped with local coordinates, one can always describe and simulate such systems locally as if they were linear spaces.
There are of course many choices of local coordinates, for rotations some famous ones are: Euler angles, the Tait-Bryan angles commonly used in aerospace applications, the unit length quaternions, and the exponentiated skew-symmetric $3\times 3$-matrices.
Lie group integrators represent a somewhat different strategy. 
Rather than specifying a choice of local coordinates from the outset, in this approach the model and the numerical integrator are
 expressed entirely in terms of a Lie group and its action on the phase space. This often leads to a more abstract and simpler formulation of the mechanical system and of the numerical schemes, deferring further details to the implementation phase. 

In the literature one can find many different types and formats of Lie group integrators. Some of these are completely general and intrinsic, meaning that they only make use of inherent properties of Lie groups and manifolds as was suggested in \cite{crouch93nio,munthe-kaas99hor, celledoni03cfl}. 
But many numerical methods have been suggested that add structure or utilise properties which are specific to a particular Lie group or manifold.
Notable examples of this are the methods based on canonical coordinates of the second kind \cite{owren00imb}, and the methods based on
the Cayley transformation \cite{lewis94caf, diele98tct}, applicable e.g. to the rotation groups and Euclidean groups.
\RE{In some applications e.g. in multi-body systems, it may be useful to formulate the problem as a mix between Lie groups and kinematic constraints, introducing for instance Lagrange multipliers. Sometimes this may lead to more practical implementations where a basic general setup involving Lie groups can be further equipped with different choices of constraints depending on the particular application. Such constrained formulations are outside the scope of the present paper. It should also be noted that the Lie group integrators devised here do not make any a priori assumptions about how the manifold is represented.}

The applications of Lie group integrators for mechanical problems also have a long history, two of the early important contributions were the Newmark methods of Simo and Vu--Quoc \cite{simo88otd} and the symplectic and energy-momentum methods by Lewis and Simo \cite{lewis94caf}. Mechanical systems are often described as Euler--Lagrange equations or as Hamiltonian systems on manifolds, with or  without  external forces, \cite{lee18gfo}. Important ideas for the discretization of mechanical systems originated also from the work of Moser and Veselov \cite{veselov88isw, moser91dvo} on discrete integrable systems. This work served as motivation for further developments in the field of geometric mechanics and for the theory of (Lie group) discrete variational integrators \cite{lee07lgv,hall17lgs,leitz18glg}. 
The majority of Lie group methods found in the literature are one-step type generalisations for classical methods, such as
Runge--Kutta type formulas. In mechanical engineering, the classical BDF methods have played an important role, and were recently
generalised \cite{wieloch19bdf} to Lie groups. Similarly,
the celebrated $\alpha$-method for linear spaces proposed by Hilber, Hughes and Taylor \cite{hilber1977ind} has been popular for solving problems in multibody dynamics, and in \RE{\cite{arnold07cot, arnold15eao, bruls10otu}} this method is generalised to a Lie group integrator. 

The literature on Lie group integrators is rich and diverse, the interested reader may consult the surveys 
 \cite{iserles00lgm,christiansen11tis, celledoni14ait, owren18lgi}
and Chapter 4 of the monograph \cite{hairer10gni} for further details.

In this paper we discuss different ways of applying Lie group integrators to simulating the dynamics of mechanical multi-body systems. %
Our point of departure is the formulation of the models as differential equations on manifolds. Assuming to be given either a Lie group acting transitively on the manifold $\mathcal{M}$ or a set of frame vector fields on $\mathcal{M}$, we use them to describe the mechanical system and further to build the numerical integrator. We shall here mostly consider schemes of the types commonly known as Crouch--Grossman methods \cite{crouch93nio}, Runge--Kutta--Munthe--Kaas methods \cite{munthe-kaas98rkm,munthe-kaas99hor} and Commutator-free Lie group methods \cite{celledoni03cfl}.

The choice of Lie group action is often not unique and thus the same mechanical system can be described in different equivalent ways. Under numerical discretization the different formulations can lead to the conservation of different geometric properties of the mechanical system. In particular, we explore the effect of these different formulations on a selection of examples in multi-body dynamics. Lie group integrators have been succesfully applied for the simulation of mechanical systems, and in problems of control, \RE{bio-mechanics} and other engineering applications, see for example \cite{park05gio}, \cite{lee07lgv} \cite{cesic16fbh}, \cite{holzinger21tio}.
The present work is 
motivated by applications in modeling and simulation of slender structures like Cosserat rods and beams \cite{simo88otd}, and one of the examples presented here is the application to a chain of pendula. Another example considers an application for the controlled dynamics of a multibody system.

In section~\ref{LGI} we give a review of the methods using only the essential intrinsic tools of Lie group integrators. The algorithms are simple and amenable for a coordinate-free description suited to object oriented implementations. In section~\ref{heavytop}, we discuss Hamiltonian systems on Lie groups, and we present three different Lie group formulations of the heavy top equations. These systems (and their Lagrangian counterpart) often arise in applications as building blocks of more realistic systems which comprise also damping and control forces. In section~\ref{implementation}, we discuss some ways of adapting the integration step size in time. In section~\ref{pendula} we consider the application to a chain of pendula.  And in section~\ref{quadrot} we consider the application of a multi-body system of interest in the simulation and control of drone dynamics.

\section{Lie group integrators} \label{LGI}

\subsection{The formulation of differential equations on manifolds}
Lie group integrators solve differential equations whose solution evolve on a manifold $\M$. \RE{For ease of notation we restrict the discussion to the case of autonomous vector fields, although allowing for explicit $t$-dependence could easily have been included.}
This means that we seek a curve $y(t)\in\M$ whose tangent at any point coincides with a vector field $F\in\mathcal{X}(\M)$ and passing through a designated initial value
$y_0$ at $t=t_0$
\begin{equation} \label{VFonM}
\dot{y}(t) = F|_{y(t)},\qquad y(t_0)=y_0.    
\end{equation}
Before addressing numerical methods for solving \eqref{VFonM} it is necessary to introduce a convenient way of representing the vector field $F$. There are different ways of doing this. 
One is to furnish $\M$ with a transitive action $\psi: G \times \M \rightarrow \M$ by some Lie group $G$ of dimension $d\geq\dim \M$. We denote the action of $g$ on $m$ as $g\cdot m$, i.e. $g\cdot m=\psi(g,m)$.
Let $\g$ be the Lie algebra of $G$, and denote by $\exp: \g\rightarrow G$ the exponential map. We define  $\infgen:\g\rightarrow\mathcal{X}(\M)$ to be the infinitesimal generator of the action, i.e.
\begin{equation} \label{frozenaction}
  \left.F_\xi\right|_m=  \left.\infgen(\xi)\right|_m = \left.\frac{d}{dt}\right|_{t=0} \psi(\exp(t\xi), m)
\end{equation}
The transitivity of the action now ensures that $\left.\infgen(\g)\right|_m=T_m\M$ for any $m\in\M$, such that any tangent vector $v_m\in T_m\M$ can be represented as $v_m=\left.\infgen(\xi_v)\right|_m$ for some $\xi_v\in\g$ ($\xi_v$ may not be unique). Consequently, for any vector field $F\in\mathcal{X}(\M)$ there exists a map $f:\M\rightarrow\g$\footnote{\RE{If the Lie group action is smooth, a map $f$ of the same regularity as $F$ can be found \cite{warner83fod}
}} such that
\begin{equation} \label{mk_generic}
    F|_m = \left.\infgen(f(m))\right|_m,\quad\text{for all}\; m\in \M
\end{equation}

This is the original tool \cite{munthe-kaas99hor} for representing a vector field on a manifold with a group action.
Another approach was used in \cite{crouch93nio} where a set of {\em frame vector fields} $E_1,\ldots, E_d$ in $\mathcal{X}(\M)$ was introduced 
assuming that for every $m\in \M$, $$\text{span}\{\left.E_1\right|_m,\ldots,\left.E_d\right|_m\}= T_m \M.$$
Then, for any vector field $F\in\mathcal{X}(\M)$ there are, in general non-unique, functions $f_i:\M\rightarrow \mathbb{R}$, 
\RE{which can be chosen with the same regularity as $F$,} such that
$$
     F|_m = \sum_{i=1}^d f_i(m) \left.E_i\right|_m.
$$
A fixed vector $\xi\in\mathbb{R}^d$ will define a vector field $F_\xi$ on $\M$ similar to \eqref{frozenaction}
\begin{equation} \label{frozenfield}
\left.F_{\xi}\right|_m = \sum_{i=1}^d \xi_i E_i|_m
\end{equation}
If $\xi_i=f_i(p)$ for some $p\in\M$, the corresponding $F_\xi$ will be a vector field in the linear span of the frame which coincides with $F$ at the point $p$. Such a vector field was named by \cite{crouch93nio} as a \emph{the vector field frozen at $p$}.

The two formulations just presented are in many cases connected, and can then be used in an equivalent manner.
Suppose that $e_1,\ldots,e_d$ is a basis of the Lie algebra $\g$, then we can simply define frame vector fields as
$E_i = \infgen(e_i)$ and the vector field we aim to describe is, 
$$
F|_m=\left.\infgen(f(m))\right|_m= \left.\infgen(\sum_i f_i(m)e_i)\right|_m=\sum_i f_i \left.E_i\right|_m.
$$
As mentioned above there is a non-uniqueness issue when defining a vector field by means of a group action or a frame.
A more fundamental description can be obtained using the machinery of connections. The assumption is that the simply connected manifold $\M$ is equipped with a connection which is flat and has constant torsion. 
Then $F_p$, the frozen vector field of $F$ at $p$ defined above, can be defined as the unique element $F_p\in\mathcal{X}(\M)$ satisfying
\begin{enumerate}
\item $F_p|_p=F|_p$
\item $\nabla_X F_p=0$ for any $X\in\mathcal{X}(M)$.
\end{enumerate}
So $F_p$ is the vector field that coincides with $F$ at $p$ and is parallel transported to any other point on $\M$ by the connection $\nabla$. Since the connection is flat, the parallel transport from the point $p$ to another point $m\in\M$ does not depend on the chosen path between the two points.
For further details, see e.g.
\cite{lundervold15oas}.
\begin{example}
For mechanical systems on Lie groups, two important constructions are the adjoint and coadjoint representatons.
For every $g\in G$ there is an automorphism $\Ad_g:\g\rightarrow\g$  defined as
$$
    \Ad_g(\xi) = TL_g\circ TR_{g^{-1}}(\xi)
$$
where $L_g$ and $R_g$ are the left and right multiplications respectively, $L_g(h)=gh$ and $R_g(h)=hg$.
Since $\Ad$ is a representation, i.e. $\Ad_{gh}=\Ad_g\circ\Ad_h$ it also defines a left Lie group action by $G$ on $\g$.
From this definition and a duality pairing $\langle\cdot,\cdot\rangle$ between $\g$ and $\g^*$, we can also derive a representation on $\g^*$ denoted $\Ad_g^*$, simply by
$$
\langle \Ad_g^*(\mu), \xi\rangle = \langle\mu, \Ad_g(\xi)\rangle,\quad \xi\in\g,\ \mu\in\g^*.
$$
The action $g\cdot \mu=\Ad_{g^{-1}}^*(\mu)$ has infinitesimal generator given as 
$$
\left.\infgen(\xi)\right|_\mu = -\ad_\xi^*\mu
$$
Following \cite{marsden94itm}, for a Hamiltonian $H:T^*G\rightarrow\mathbb{R}$, define $H^-$ to be its restriction to $\g^*$.
Then the Lie-Poisson reduction of the dynamical system is defined on $\g^*$ as
$$
\RE{\dot{\mu} = -\ad^*_{\frac{\partial H^-}{\partial \mu}}\mu}
$$
and this vector field is precisely of the form \eqref{mk_generic} with \RE{$f(\mu)=\frac{\partial H^-}{\partial\mu}(\mu)$}.
A side effect of this is that the integral curves of these Lie-Poisson systems preserve coadjoint orbits, making the
coadjoint action an attractive choice for Lie group integrators.

Let us now detail the situation for the very simple case where $G=SO(3)$. The Lie algebra $\mathfrak{so}(3)$ can be
modeled as $3\times 3$ skew-symmetric matrices, and via the standard basis we identify each such matrix $\hat{\xi}$ by
a vector $\xi\in\mathbb{R}^3$, this identification is known as the hat map
\begin{equation} \label{hatmap}
\hat{\xi} = \left[
\begin{array}{ccc}
  0   & -\xi_3 & \xi_2 \\
  \xi_3   & 0 & -\xi_1 \\
  -\xi_2 & \xi_1 & 0
\end{array}
\right]
\end{equation}
Now, we also write the elements of $\mathfrak{so}(3)^*$ as vectors in $\mathbb{R}^3$ with duality pairing $\langle\mu,\xi\rangle=\mu^T\xi$. With these representations, we find that the coadjoint action can be expressed as
$$
g\cdot\mu=\psi(g,\mu) = \Ad_{g^{-1}}^*\mu = g\mu
$$
the rightmost expression being a simple matrix-vector multiplication. Since $g$ is orthogonal, it follows that the coadjoint orbits foliate 3-space into spherical shells, and the coadjoint action is transitive on each of these orbits.
The free rigid body can be cast as a problem on $TSO(3)^*$ with a left invariant Hamiltonian which reduces to the
function 
$$
H^-(\mu)=\frac12 \langle \mu, \mathbb{I}^{-1}\mu\rangle
$$
on $\mathfrak{so}(3)^*$ where $\mathbb{I}:\mathfrak{so}(3)\rightarrow\mathfrak{so}(3)^*$ is the inertia tensor.
From this, we can now set \RE{$f(\mu)=\partial H^-/\partial \mu = \mathbb{I}^{-1}\mu$}. \RE{We then} recover the Euler free rigid body equation as
$$
\dot{\mu} = \left.\infgen(f(\mu)\right|_\mu = -\ad_{\mathbb{I}^{-1}\mu}^*\mu = -\mathbb{I}^{-1}\mu \times \mu
$$
where the last expression involves the cross product of vectors in $\mathbb{R}^3$.
\end{example}

\subsection{Two classes of Lie group integrators} \label{twoclasses}
The simplest numerical integrator for linear spaces is the explicit Euler method. Given an initial value problem $\dot{y}=F(y)$, $y(0)=y_0$ the method is defined as $y_{n+1}=y_n + hF(y_n)$ for some stepsize $h$. In the spirit of the previous section, one could think of
the Euler method as the $h$-flow of the constant vector field $F_{y_n}(y)=F(y_n)$, that is
$$
y_{n+1} = \exp(hF_{y_n})\,y_n
$$
This definition of the Euler method makes sense also when $F$ is replaced by a vector field on some manifold. In this general situation it is known as the Lie--Euler method.

We shall here consider the two classes of methods known as
Runge--Kutta--Munthe--Kaas (RKMK) methods and Commutator-free Lie group methods.

For RKMK methods the underlying idea is to transform the problem from the manifold $\M$ to the Lie algebra $\g$, take a time step, and map the result back to $\M$. The transformation we use is
$$
y(t) = \exp(\sigma(t))\cdot y_0,\quad\sigma(0)=0.
$$
The transformed differential equation for $\sigma(t)$ makes use of the derivative of the exponential mapping, the reader should consult
\cite{munthe-kaas99hor} for details about the derivation, we give the final result
\begin{equation}\label{dexpinveq}
\dot{\sigma}(t) = \dexp_{\sigma(t)}^{-1} (f(\exp(\sigma(t))\cdot y_0))
\end{equation}
The map $v\mapsto\dexp_u(v)$ is linear and invertible when $u$ belongs to some sufficiently small neighborhood of $0\in\g$. 
It has an expansion in nested Lie brackets \cite{hausdorff06dse}. Using the operator $\ad_u(v)=[u,v]$ and its powers
$\ad_u^2 v=[u,[u,v]]$ etc, one can write
\begin{equation} \label{dexpseries}
\dexp_u(v) = \left.\frac{e^z-1}{z}\right|_{z=\ad_u}(v) = v + \frac12[u,v] + \frac16[u,[u,v]] + \cdots
\end{equation}
and the inverse is
\begin{equation} \label{dexpinvseries}
 \dexp_u^{-1}(v) =\left.\frac{z}{e^z-1}\right|_{z=\ad_u}(v)= v -\frac12[u,v] + \frac1{12}[u,[u,v]]+\cdots
\end{equation}
The RKMK methods are now obtained simply by applying some standard Runge--Kutta method to the transformed equation \eqref{dexpinveq} with a time step $h$, using initial value $\sigma(0)=0$. This leads to an output $\sigma_1\in\g$ and one simply sets $y_1=\exp(\sigma_1)\cdot y_0$. Then one repeats the procedure replacing $y_0$ by $y_1$ in the next step etc. While solving \eqref{dexpinveq} one needs to evaluate $\dexp_u^{-1}(v)$ as a part of the process. This can be done by truncating the series \eqref{dexpinvseries} since $\sigma(0)=0$ implies that we always evaluate $\dexp_u^{-1}$ with $u=\mathcal{O}(h)$, and thus, the $k$th iterated commutator $\ad_u^k=\mathcal{O}(h^k)$.
For a given Runge--Kutta method, there are some clever tricks that can be done to minimise the total number of commutators to be included from the expansion of $\dexp_u^{-1}v$, see \cite{casas03cel, munthe-kaas99cia}. We give here one concrete example of an RKMK method proposed in \cite{casas03cel}
\begin{align*}
   f_{n,1} &= h f(y_n),\\
   f_{n,2} &= h f(\exp(\tfrac{1}{2}f_{n,1}) \cdot y_n), \\
   f_{n,3} &= h f(\exp(\tfrac{1}{2}f_{n,2}-\tfrac{1}{8}[f_{n,1},f_{n,2}])\cdot y_n), \\
   f_{n,4} &= h f(\exp(f_{n,3})\cdot y_n), &  \\
   y_{n+1} &= \exp(\tfrac{1}{6}(f_{n,1}+2f_{n,2}+2f_{n,3}+f_{n,4}-\tfrac12[f_{n,1},f_{n,4}]))\cdot y_n.
\end{align*}

The other option is to compute the exact expression for $\dexp_u^{-1}(v)$ for the particular Lie algebra we use. For instance, it was shown in \cite{celledoni03lgm} that for the Lie algebra $\so(3)$ one has
$$
\dexp_u^{-1}(v)=v - \frac12 u\times v + \alpha^{-2}(1-\tfrac{\alpha}{2}\cot\tfrac{\alpha}{2})\; u\times (u\times v)
$$
We will present the corresponding formula for $\se(3)$ in Section~\ref{dexpinvse3}.

The second class of Lie group integrators to be considered here are the commutator-free methods, named this way in \cite{celledoni03cfl} to emphasize the contrast to RKMK schemes which usually include commutators in the method format. These schemes include the Crouch-Grossman methods \cite{crouch93nio} and they have the format
\begin{align*}
    Y_{n,r} &= \exp\left(h\sum_{k}\alpha_{r,J}^k f_{n,k}\right)\cdots \exp\left(h\sum_{k}\alpha_{r,1}^k f_{n,k}\right)\cdot y_n \\
      f_{n,r} &= f(Y_{n,r}) \\[1mm]
    y_{n+1} &= \exp\left(h\sum_k \beta_J^k f_{n,k}\right)\cdots \exp\left(h\sum_k \beta_1^k f_{n,k}\right)\cdot y_n
\end{align*}
Here the Runge--Kutta coefficients $\alpha_{r,j}^k$, $\beta_{j}^r$ are related to a classical Runge--Kutta scheme with coefficients $a_r^k$, $b_r$ in that $a_r^k=\sum_j \alpha_{r,j}^k$ and $b_r=\sum_j \beta_{j}^r$. The $\alpha_{r,j}^k$, $\beta_{j}^r$ are usually chosen to obtain computationally inexpensive schemes with the highest possible order of convergence. The computational complexity of the above schemes depends on the cost of computing an exponential as well as of evaluating the vector field. Therefore it makes sense to keep the number of exponentials $J$ in each stage as low as possible, and possibly also the number of stages $s$. A trick proposed in \cite{celledoni03cfl} was to select coefficients \RE{that make it possible to reuse exponentials} from one stage to another. This is perhaps best illustrated through the following example from \cite{celledoni03cfl}, a generalisation of the classical 4th order Runge--Kutta method.
\begin{align}
\begin{split} \label{cf4}
Y_{n,1} &= y_n \\
Y_{n,2} &=\exp(\tfrac12 hf_{n,1}) \cdot y_n \\
Y_{n,3} &= \exp(\tfrac12 hf_{n,2}) \cdot y_n \\
Y_{n,4} &= \exp(h f_{n,3}-\tfrac12 h f_{n,1}) \cdot Y_{n,2} \\
y_{n+\frac12} &=\exp(\tfrac{1}{12}h(3f_{n,1}+2f_{n,2}+2f_{n,3}-f_{n,4})) \cdot y_n \\
y_{n+1} &=\exp(\tfrac{1}{12}h(-f_{n,1}+2f_{n,2}+2f_{n,3}+ 3f_{n,4})) \cdot y_{n+\frac12}
\end{split}
\end{align}
where $f_{n,i}=f(Y_{n,i})$. Here, we see that one exponential is saved in computing $Y_{n,4}$ by making use of $Y_{n,2}$.

\subsection{An exact expression for $\dexp_u^{-1}(v)$ in $\mathfrak{se}(3)$} \label{dexpinvse3}

As an alternative to using a truncated version of the infinite series for $\dexp_u^{-1}$ \eqref{dexpinvseries}, one can consider exact expressions obtained for certain Lie algebras. Since $\mathfrak{se}(3)$ is particularly important in applications to mechanics, we give here its exact expression. For this, we represent elements of $\mathfrak{se}(3)$ as a
pair $(A,a)\in\mathbb{R}^3\times\mathbb{R}^3\cong\mathbb{R}^6$, the first component corresponding to a skew-symmetric matrix $\hat{A}$
via \eqref{hatmap} and $a$ is the translational part. Now, let $\varphi(z)$ be a real analytic function at $z=0$. We define
$$
\varphi_+(z) = \frac{\varphi(iz)+\varphi(-iz)}{2},\qquad \varphi_-(z) = \frac{\varphi(iz)-\varphi(-iz)}{2i}
$$
We next define the four functions
$$
g_1(z) = \frac{\varphi_-(z)}{z},\ \tilde{g}_1(z) = \frac{g_1'(z)}{z},\quad 
g_2(z) = \frac{\varphi(0)-\varphi_+(z)}{z^2},\ \tilde{g}_2(z)=\frac{g_2'(z)}{z} 
$$
and the two scalars $\rho=A^Ta$, $\alpha=\|A\|_2$. One can show that for any $(A,a)$ and $(B,b)$ in $\mathfrak{se}(3)$, it holds that
$$
\varphi(\ad_{(A,a)})(B,b) = (C,c)
$$
where
\begin{align*}
C&=\varphi(0)B + g_1(\alpha) A\times B + g_2(\alpha)\, A\times (A\times B)\\
c&=\varphi(0)b + g_1(\alpha)\, (a\times B+A\times b)
+\rho\tilde{g}_1(\alpha)\,A\times B  + \rho\tilde{g}_2(\alpha)\, A\times (A\times B)\\
&+ g_2(\alpha)\, (a\times (A\times B)+A\times (a\times B) + A\times (A\times b))
\end{align*} 
Considering for instance \eqref{dexpinvseries}, we may now use $\varphi(z)=\frac{z}{e^z-1}$ to calculate
$$
g_1(z) = -\frac12,\ \tilde{g}_1(z)=0,\ g_2(z) = \frac{1-\tfrac{z}{2}\cot\tfrac{z}{2}}{z^2},\ 
\tilde{g}_2(z)=\frac1z\frac{d}{dz}g_2(z),\ \varphi(0)=1.
$$
and thereby obtain an expression for $\dexp_{(A,a)}^{-1}(B,b)$ with the formula above.

Similar types of formulas are known for computing the matrix exponential as well as functions of the $\ad$-operator for several other Lie groups of small and medium dimension. For instance in \cite{mueller17cmf} a variety of coordinate mappings for rigid body motions are discussed.
For Lie algebras of larger dimension, both the exponential mapping and  
$\dexp_u^{-1}$ may become computationally infeasible. For these cases, one may benefit from replacing the exponential by some other coordinate map for the Lie group $\phi:\g\rightarrow G$. One option is to use canonical coordinates of the second kind \cite{owren00imb}. Then for some Lie groups such as the orthogonal, unitary and symplectic groups, there exist other maps that can be used and which are computationally less expensive. A popular choice is the Cayley transformation \cite{diele98tct}.

\section{Hamiltonian systems on Lie groups}~\label{heavytop}

In this section we consider Hamiltonian systems on Lie groups. These systems (and their Lagrangian counterpart) often appear in mechanics applications as building blocks for more realistic systems with additional damping and control forces. We consider canonical systems on the cotangent bundle of a Lie group and Lie-Poisson systems which can arise by symmetry reduction or otherwise. We illustrate the various cases with different formulations of the heavy top system.

\subsection{Semi-direct products } 

The coadjoint action by $G$ on $\mathfrak{g}^*$ is denoted 
$\mathrm{Ad}_g^*$ defined for any $g\in G$ as
\begin{equation} \label{eq:coadjoint-action}
\langle \mathrm{Ad}_g^*\mu, \xi\rangle = \langle \mu, \mathrm{Ad}_g\xi\rangle,\quad\forall\xi\in\mathfrak{g},
\end{equation}
where  $\mathrm{Ad}:\mathfrak{g}\rightarrow\mathfrak{g}$ is the adjoint representation and
for a duality pairing $\langle {\cdot}, {\cdot} \rangle$ between $\mathfrak{g}^*$ and $\mathfrak{g}$.

We consider the cotangent bundle of a Lie group $G$, $T^*G$ and identify it with  $G\times \mathfrak{g}^*$ using the right multiplication $R_g:G\rightarrow G$ and its tangent mapping $R_{g*}:=TR_g$. The cartesian product $G\times \mathfrak{g}^*$ can be given a semi-direct product structure that turns it into a Lie group \RE{$\mathbf{G} :=G\ltimes\mathfrak{g}^*$} where the group multiplication is
\begin{equation}
(g_1,\mu_1)\cdot (g_2,\mu_2)=(g_1\cdot g_2, \mu_1+\mathrm{Ad}^*_{g_1^{-1}}\mu_2).
\end{equation}

Acting by left multiplication any vector field $F\in\mathcal{X}(\mathbf{G})$ is expressed by means of a map $f\colon \mathbf{G} \rightarrow T_e \mathbf{G}$,
\begin{equation} \label{eq:Fpres}
        F(g,\mu) = T_e R_{(g,\mu)} f(g,\mu) = (R_{g*} f_1, f_2-\mathrm{ad}_{f_1}^*\mu),
\end{equation}
where $f_1=f_1(g,\mu)\in\mathfrak{g}$, $f_2=f_2(g,\mu)\in\mathfrak{g}^*$ are the two components of $f$.

\subsection{Symplectic form and Hamiltonian vector fields}

The right trivialised\footnote{$ \omega_{(g,\mu)}$ is obtained from
the natural symplectic form on
$T^*G$ (which is a differential two-form), defined as
\[
\Omega_{(g,p_g)}((\delta v_1, \delta \pi_1),(\delta v_2,\delta\pi_2))
=\langle \delta \pi_2, \delta v_1\rangle - \langle \delta \pi_1, \delta v_2\rangle,
\]
by right trivialization.
} symplectic form pulled back to $\mathbf{G}$ reads
\RE{{
\begin{equation} \label{eq:sympform}
     \begin{split}
      \omega_{(g,\mu)}( (R_{g*} \xi_1, \delta\nu_1), (R_{g*}\xi_2, \delta\nu_2))
      &= \langle\delta\nu_2,\xi_1 \rangle  + \\
      &-\langle\delta\nu_1, \xi_2 \rangle - \langle\mu, [\xi_1,\xi_2]\rangle, \qquad \xi_1,\xi_2 \in \mathfrak{g}.
      \end{split}
\end{equation}
}}
See \cite{lewis94caf} for more details, proofs and for a the left trivialized symplectic form.
The vector field $F$ is a Hamiltonian vector field if it satisfies
$$\mathrm{i}_F\omega=dH,$$
for some Hamiltonian function $H:T^*G\rightarrow \mathbb{R}$, where $\mathrm{i}_F$ is defined as $\mathrm{i}_F(X):=\omega(F,X)$ for any vector field $X$. This implies that the map $f$ for such a Hamiltonian vector field gets the form
\begin{equation}
\label{hvfongstar}
f(g, \mu)=\left(  \frac{\partial H}{\partial \mu}(g,\mu), -R_g^*\frac{\partial H}{\partial g}(g,\mu)\right).
\end{equation}

The following is a one-parameter family of symplectic Lie group integrators on $T^*G$:

\begin{align}
\label{symplectic0}
M_{\theta}&=\mathrm{dexp}_{-\xi}^{*} (\mu_0 + \mathrm{Ad}^*_{\exp(\theta \xi)}(\bar{n})) - \theta \mathrm{dexp}_{-\theta\xi}^{*}  \mathrm{Ad}^*_{\exp(\theta \xi)}(\bar{n}),\\
\label{symplectic1}
	(\xi, \bar n) &= h f\bigl(\exp(\theta \xi) \cdot g_0, M_{\theta} \bigr), \\
	\label{symplectic2}
	(g_1, \mu_1) &= (\exp(\xi), \mathrm{Ad}_{\exp((\theta-1) \xi)}^{*} \bar n) \cdot (g_0, \mu_0).
\end{align}

For higher order integrators of this type and a complete treatment see \cite{bogfjellmo15hos}.

\subsection{Reduced equations Lie Poisson systems}

A mechanical system formulated on the cotangent bundle $T^*G$
with a left or right invariant Hamiltonian can be reduced to a system
on $\mathfrak{g}^*$ \cite{marsden84sdp}. In fact for a Hamiltonian $H$ right invariant under the left action of $G$,  $\frac{\partial H}{\partial g}=0$, and from \eqref{eq:Fpres} and \eqref{hvfongstar} we get for the second equation
\begin{equation}
\label{Lie-Poisson}
     \dot{\mu} = \mp\mathrm{ad}^*_{\frac{\partial H}{\partial\mu}} \mu,
\end{equation}
where the positive sign is used in case of left invariance
(see e.g.\ section~13.4 in \cite{marsden99itm}). 
The solution to this system preserves coadjoint orbits, thus using the Lie group action 
\[
    g\cdot\mu = \mathrm{Ad}_{g^{-1}}^*\mu,
\]
to build a Lie group integrator results in preservation of such coadjoint orbits. 
Lie group integrators for this interesting case were studied in \cite{engo01nio}.

The Lagrangian counterpart to these Hamiltonian equations are \RE{the} Euler--Poincar{\'e} equations\footnote{The Euler--Poincar{\'e} equations are Euler--Lagrange equations with respect to a Lagrange--d'Alembert principle obtained taking constraint variations.}, \cite{holm98tep}.


\subsection{Three different formulations of the heavy top equations}

The heavy top is a simple test example for illustrating the behaviour of Lie group methods. We will consider three different formulations for this mechanical system. The first formulation is on  $T^*SO(3)$ where the equations are canonical Hamiltonian, a second point of view is that the system is a Lie--Poisson system on $\mathfrak{se}(3)^*$, 
and finally it is canonical Hamiltonian on a larger group with a quadratic Hamiltonian function. 
The three different formulations suggest the use of different Lie group integrators.

\begin{figure}[h]
    \centering
    \includegraphics[width=.6\textwidth]{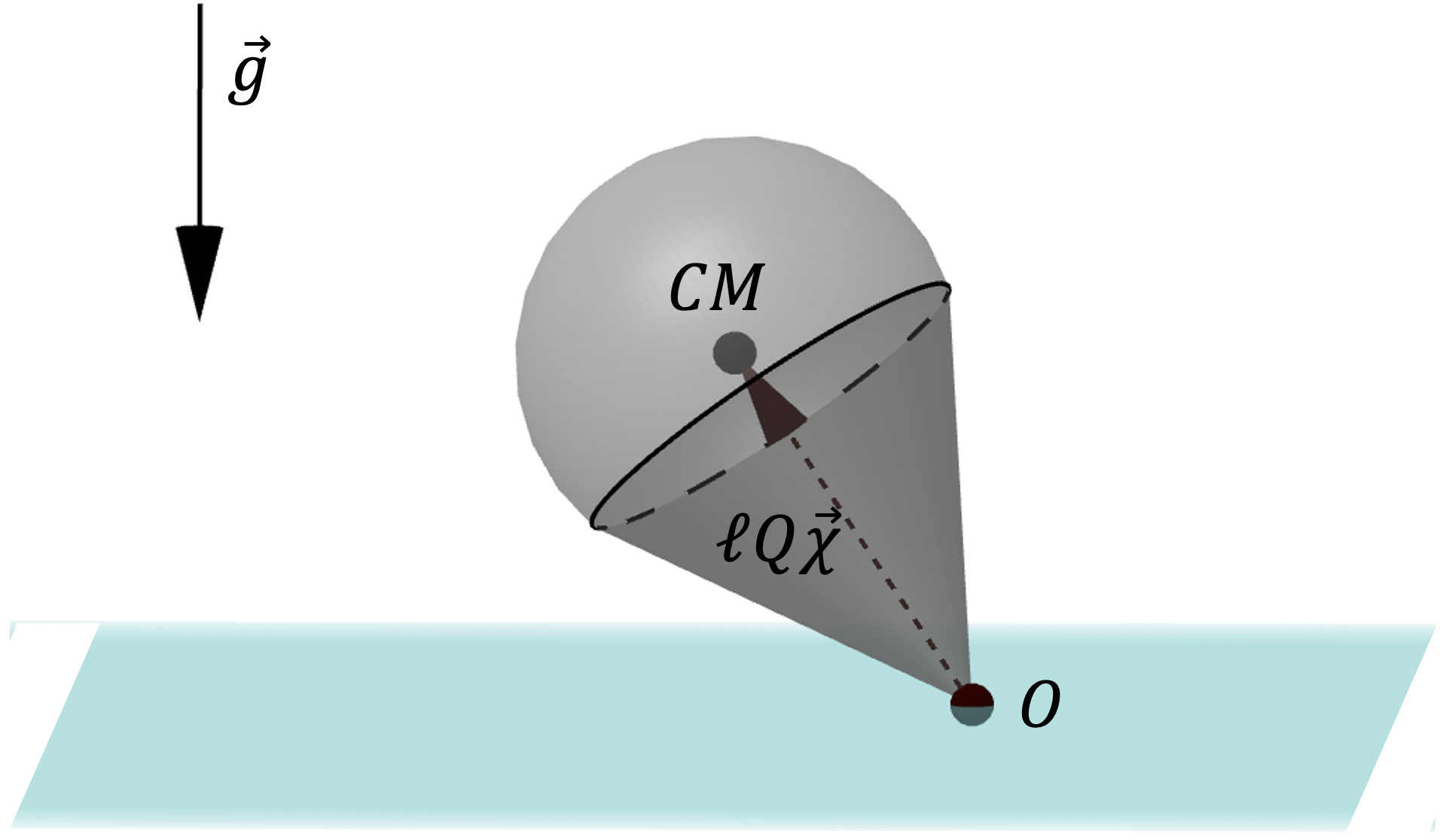}
    \caption{\RE{Illustration of the heavy top, where $CM$ is the center of mass of the body, $O$ is the fixed point, $\vec{g}$ is the gravitational acceleration vector, and $\ell,Q,\vec{\chi}$ follow the notation introduced in Section \ref{heavytop1}}}
\end{figure}

\subsubsection{Heavy top equations on $T^*SO(3)$.}\label{heavytop1}

The heavy top is a rigid body with a fixed point in a gravitational field. The phase space of this mechanical system is $T^*SO(3)$ where the equations of the heavy top are in canonical Hamiltonian form. Assuming $(Q,p)$ are coordinates for $T^*SO(3)$, $\Pi=(T_eL_Q)^*(p)$ is the left trivialized or body momentum. 
The Hamiltonian of the heavy top is given in terms of $(Q,\Pi)$  as
\[
    H \colon SO(3)\ltimes\mathfrak{so}(3)^* \to \mathbb{R}, \quad H(Q,\Pi) =\frac{1}{2} \langle \Pi, \mathbb{I}^{-1}\Pi\rangle + Mg\ell \,\mathbf{\Gamma} \cdot \mathcal{X},\quad \mathbf{\Gamma}=Q^{-1}\mathbf{\Gamma}_0,
\]
where $\mathbb{I}\colon \mathfrak{so}(3)\rightarrow \mathfrak{so}(3)^*$ is the inertia tensor, here represented as a diagonal \label{vale:pgref}
$3\times 3$ matrix, $\mathbf{\Gamma}=Q^{-1}\mathbf{\Gamma}_0$, where $\mathbf{\Gamma}_0\in \mathbb{R}^3$ 
is the axis of the spatial coordinate system parallel to the direction of gravity but pointing upwards,
$M$ is the mass of the body, $g$ is the gravitational acceleration, $\mathcal{X}$ is  the body fixed unit vector of the oriented line segment pointing from the fixed point to the center of mass of the body, $\ell$ is the length of this segment. The equations of motion on $SO(3)\ltimes\mathfrak{so}(3)^*$ are
\begin{align}
\label{ht1}
\dot{\Pi}&=\Pi \times \mathbb{I}^{-1}\Pi + Mg\ell \, \mathbf{\Gamma}\times \mathcal{X},\\
\label{ht2}
\dot{Q} &=Q\, \widehat{ \mathbb{I}^{-1}\Pi }.
\end{align}

The identification of $T^*SO(3)$ with $SO(3)\ltimes\mathfrak{so}(3)^*$ via right trivialization leads to the spatial \RE{momentum} variable $\pi=(T_eR_Q)^*(p)=Q\Pi$. The equations written in the space variables $(Q,\pi)$ get the form

\begin{align}
\label{ht10}
\dot{\pi}&= Mg\ell \, \mathbf{\Gamma}_0\times Q\mathcal{X},\\
\label{ht20}
\dot{Q} &=\hat{\omega} Q \qquad \omega =Q \mathbb{I}^{-1}Q^T\pi .
\end{align}
where, the first equation states that the component of $\pi$ parallel to $\mathbf{\Gamma}_0$ is constant in time.
These equations can be obtained from \eqref{eq:Fpres} and \eqref{hvfongstar} on the right trivialized $T^*SO(3)$, $SO(3)\ltimes\mathfrak{so}(3)^*$,  with the heavy top Hamiltonian and the symplectic Lie group integrators \eqref{symplectic1}-\eqref{symplectic2} can be applied in this case. Similar methods were proposed in \cite{lewis94caf} and \cite{saccon09mrf}.


\subsubsection{Heavy top equations on $\mathfrak{se}^*(3)$}

The Hamiltonian of the heavy top is not invariant under the action of $SO(3)$, so the equations \RE{\eqref{ht1}-\eqref{ht2} given} in section~\eqref{heavytop1} cannot be reduced to $\mathfrak{so}^*(3)$, nevertheless the heavy top equations are Lie--Poisson on $\mathfrak{se}^*(3)$, \cite{vinogradov77tso,guillemin80tmm,ratiu81epe}.

Observe that the equations of the heavy top on $T^*SO(3)$ \eqref{ht1}-\eqref{ht2} can be easily modified eliminating the variable $Q\in SO(3)$ and replacing it with $\mathbf{\Gamma}\in \mathbb{R}^3$ $\mathbf{\Gamma}=Q^{-1}\mathbf{\Gamma}_0$ to obtain
\begin{align}
\label{ht11}
\dot{\Pi}&=\Pi \times \mathbb{I}^{-1}\Pi + Mg\ell \, \mathbf{\Gamma}\times \mathcal{X},\\
\label{ht12}
\dot{\mathbf{\Gamma}} &=\mathbf{\Gamma}\times  (\mathbb{I}^{-1}\Pi ).
\end{align}
We will see that the solutions of these equations evolve on $\mathfrak{se}^*(3)$.
In what follows, we consider elements of $\mathfrak{se}^*(3)$ to be pairs of vectors in $\mathbb{R}^3$, e.g. $(\Pi, \mathbf{\Gamma})$. Correspondingly the elements of $SE(3)$ are represented as pairs $(g,\mathbf{u})$ with $g\in SO(3)$ and $\mathbf{u}\in \mathbb{R}^3$.
The group multiplication in  $SE(3)$ is then
$$(g_1,\mathbf{u}_1)\cdot (g_2,\mathbf{u}_2)=\RE{(g_1g_2,g_1\mathbf{u}_2 + \mathbf{u}_1)},$$
where $g_1g_2$ is the product in $SO(3)$ and $g_1\mathbf{u}$ is the product of a $3\times 3$ orthogonal matrix with a vector in $\mathbb{R}^3$.
The coadjoint representation and its infinitesimal generator on $\mathfrak{se}^*(3)$ take the form
$$\mathrm{Ad}^*_{(g,\mathbf{u})}(\Pi, \mathbf{\Gamma})=(g^{-1}(\Pi-\mathbf{u}\times \mathbf{\Gamma}), g^{-1}\mathbf{\Gamma}),\;\; \mathrm{ad}^*_{(\xi, \mathbf{u})}(\Pi, \mathbf{\Gamma})=(-\xi \times \Pi-\mathbf{u}\times \mathbf{\Gamma}, -\xi \times \mathbf{\Gamma}).$$
Using this expression for $\mathrm{ad}^*_{(\xi, \mathbf{u})}$ with $(\xi=\frac{\partial H}{\partial \Pi}, \mathbf{u}=\frac{\partial H}{\partial  \mathbf{\Gamma}})$, it can be easily seen that the equations \eqref{Lie-Poisson} in this setting reproduce the heavy top equations \eqref{ht11}-\eqref{ht12}.
Therefore the equations are Lie--Poisson equations on $\mathfrak{se}^*(3)$. 
However since the heavy top is a rigid body with a fixed point and there are no translations, these equations do not arise from a reduction of $T^*SE(3)$. 
Moreover the Hamiltonian on $\mathfrak{se}(3)^*$ is not quadratic and the equations are not geodesic equations. 
Implicit and explicit Lie group integrators applicable to this formulation of the heavy top equations and preserving coadjoint orbits were discussed in \cite{engo01nio}, for a variable stepsize integrator applied to this formulation of the heavy top see \cite{curry19vss}.


\subsubsection{Heavy top equations with quadratic Hamiltonian.}

We rewrite the heavy top equations one more time considering the constant vector $\mathbf{p}=-Mg\ell \mathcal{X}$ as a momentum variable conjugate to the position $\mathbf{q}\in \mathbb{R}^3$ and where $\mathbf{p}=Q^{-1}\mathbf{\Gamma}_0+\dot{\mathbf{q}}$, and the Hamiltonian is a quadratic function of $\Pi$, $Q$, $\mathbf{p}$ and $\mathbf{q}$:
\begin{align*}
&H \colon T^*SO(3)\times {\mathbb{R}^3}^*\times {\mathbb{R}^3} \to \mathbb{R}, \\ 
&H((\Pi,Q), (\mathbf{p}, \mathbf{q})) =\frac{1}{2} \langle \Pi, \mathbb{I}^{-1}\Pi\rangle + \frac{1}{2}\|\mathbf{p}- Q^{-1}\mathbf{\Gamma}_0\|^2-\frac{1}{2}\|Q^{-1}\mathbf{\Gamma}_0\|^2,
\end{align*}
see \cite[section 8.5]{Holm}.
This Hamiltonian is invariant under the left action of $SO(3)$. The corresponding equations are canonical on $T^*S\equiv S\ltimes \mathfrak{s}^*$ where $S=SO(3)\times \mathbb{R}^3$ with Lie algebra $\mathfrak{s}:=\mathfrak{so}(3)\times {\mathbb{R}^3}$ and $T^*S$ can be identified with $T^*SO(3)\times {\mathbb{R}^3}^* \times \mathbb{R}^3$. The equations are
\begin{align}
\label{ht21}
\dot{\Pi}&=\Pi \times \mathbb{I}^{-1}\Pi -  (Q^{-1}\mathbf{\Gamma}_0) \times \mathbf{p},\\
\label{ht22}
\dot{Q} &=Q\, \widehat{ \mathbb{I}^{-1}\Pi },\\
\label{ht23}
\dot{\mathbf{p}}&=\RE{\mathbf{0}},\\
\dot{\mathbf{q}}&= \mathbf{p}-Q^{-1}\mathbf{\Gamma}_0.
\end{align}
and in the spatial momentum variables
\begin{align}
\label{ht211}
\dot{\pi}&= -  \mathbf{\Gamma}_0 \times Q\mathbf{p},\\
\label{ht221}
\dot{Q} &=\hat{\omega}Q,\quad  \omega
=Q \mathbb{I}^{-1}Q^T\pi ,\\
\label{ht231}
\dot{\mathbf{p}}&=\RE{\mathbf{0}},\\
\label{ht241}
\dot{\mathbf{q}}&= \mathbf{p}-Q^{-1}\mathbf{\Gamma}_0.
\end{align}
Similar formulations were considered in \cite{leonard97sad} for the stability analysis of an underwater vehicle. A similar but different formulation of the heavy top was considered in \cite{bruls10otu}.

\subsubsection{Numerical experiments.}

We apply various implicit Lie group integrators to the heavy top system. The test problem we consider is the same as in \cite{bruls10otu}, where $Q(0)=I$, $\ell=2$, $M=15$ \newline $\mathbb{I}=\mathrm{diag}(0.234375, 0.46875, 0.234375)$, $\pi(0)=\mathbb{I}(0, 150, -4.61538)$, $\mathcal{X}=(0,1,0)$ $\Gamma_0=(0,0,-9.81)$.

\begin{figure}[htbp]
    \centering
    \includegraphics[width=0.31\textwidth]{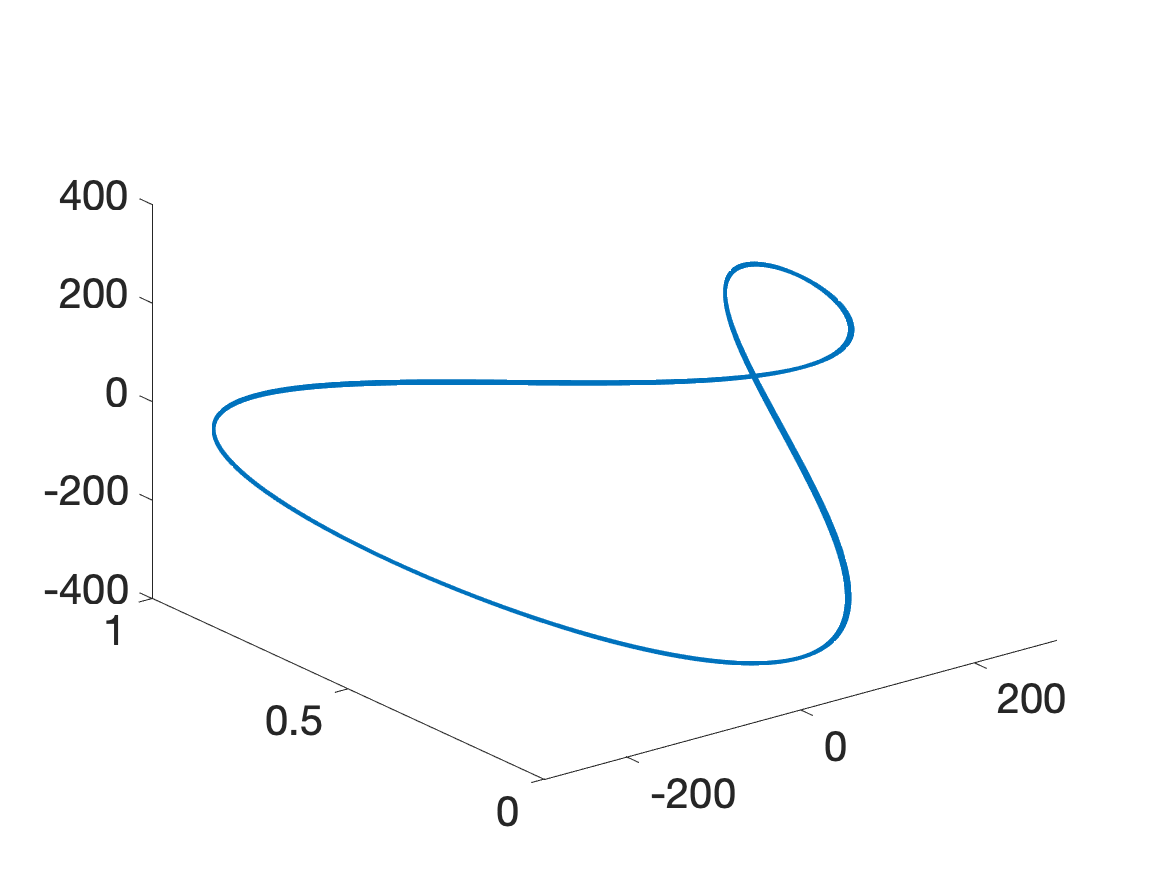}
    \includegraphics[width=0.31\textwidth]{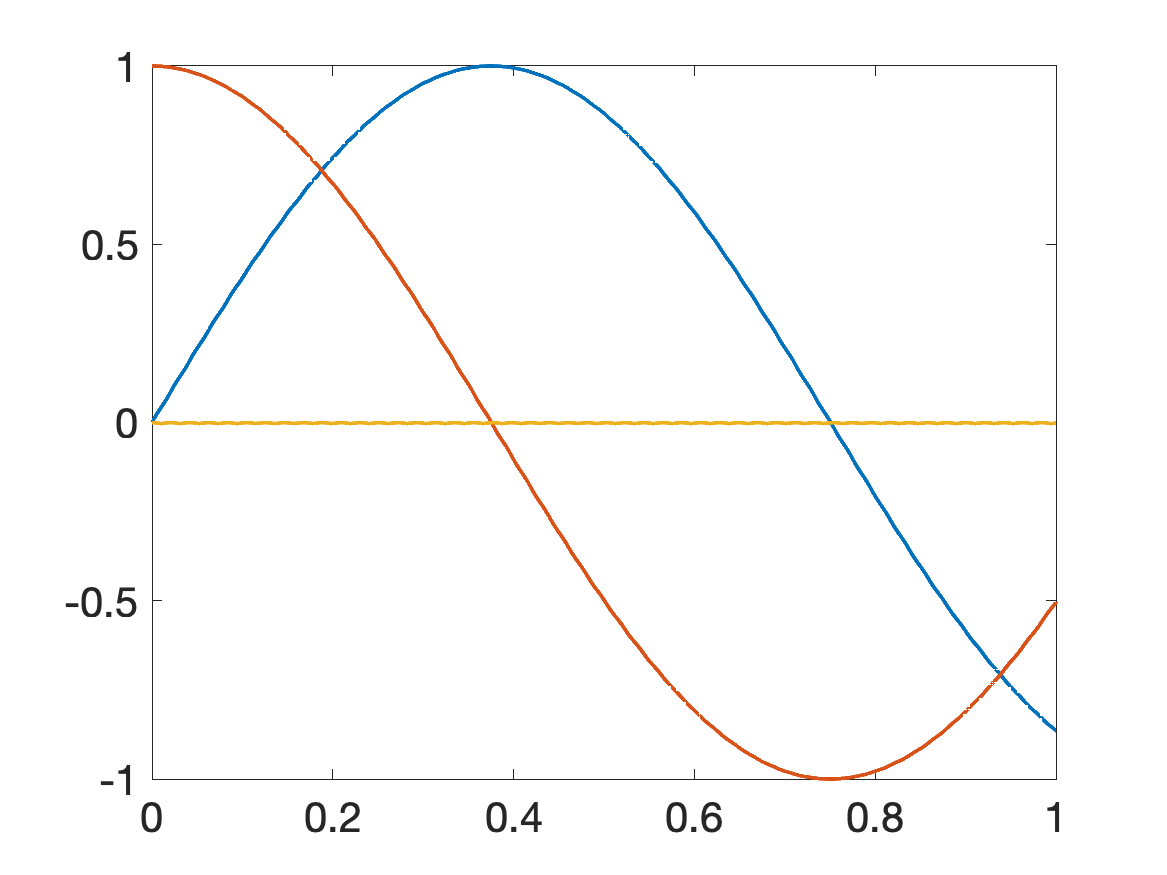}
    \includegraphics[width=0.35\textwidth]{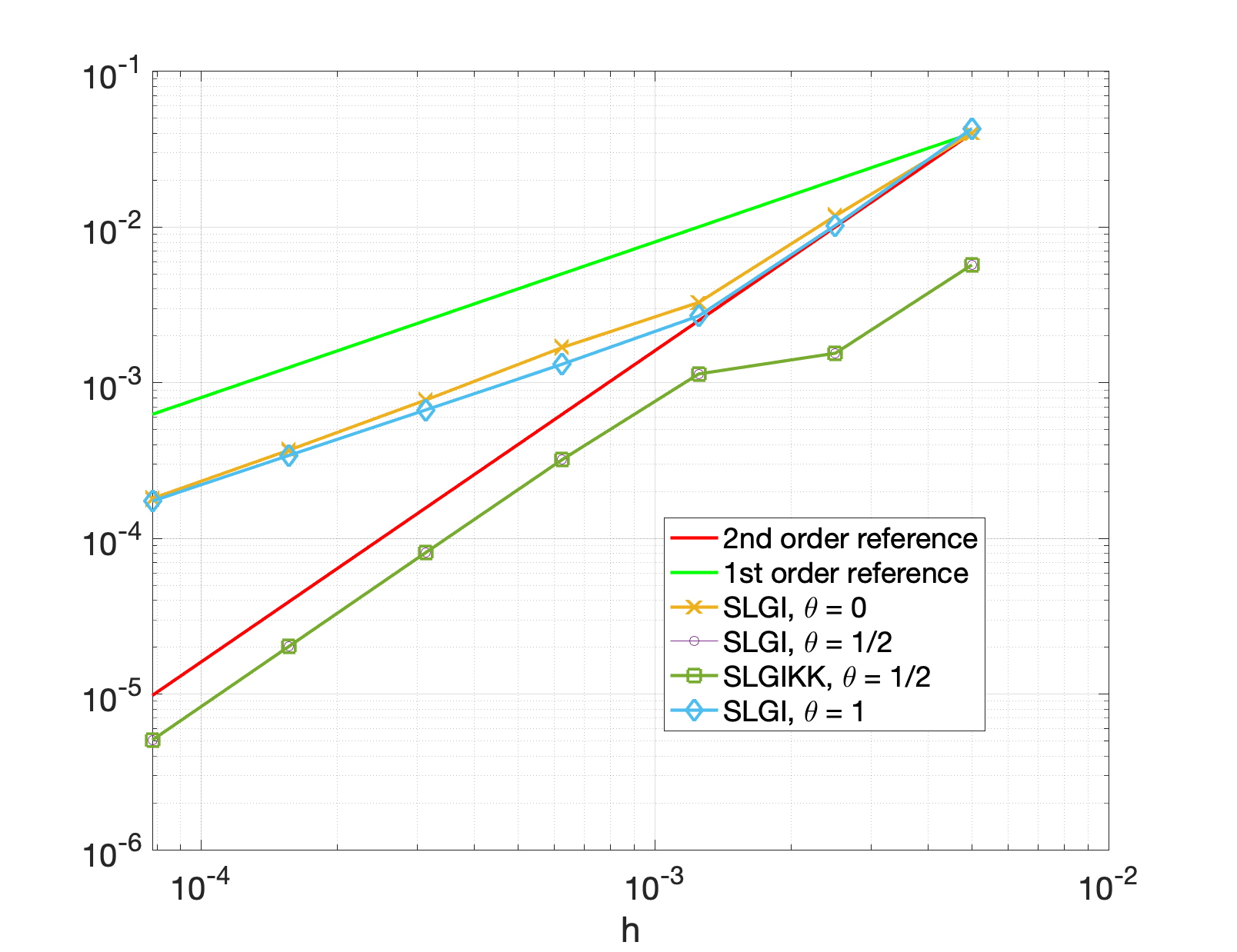}
    \caption{Symplectic Lie group integrators integration on the time interval $[0,1]$. Left:  3D plot of $M\ell Q^{-1}\Gamma_0$. Center: components of $Q\mathcal{X}$. \RE{The left and center plots are computed with the same  step-size.} Right: verification of the order of the methods.}
    \label{fig:HT1}
\end{figure}
\begin{figure}[htbp]
    \centering
    \includegraphics[width=0.9\textwidth]{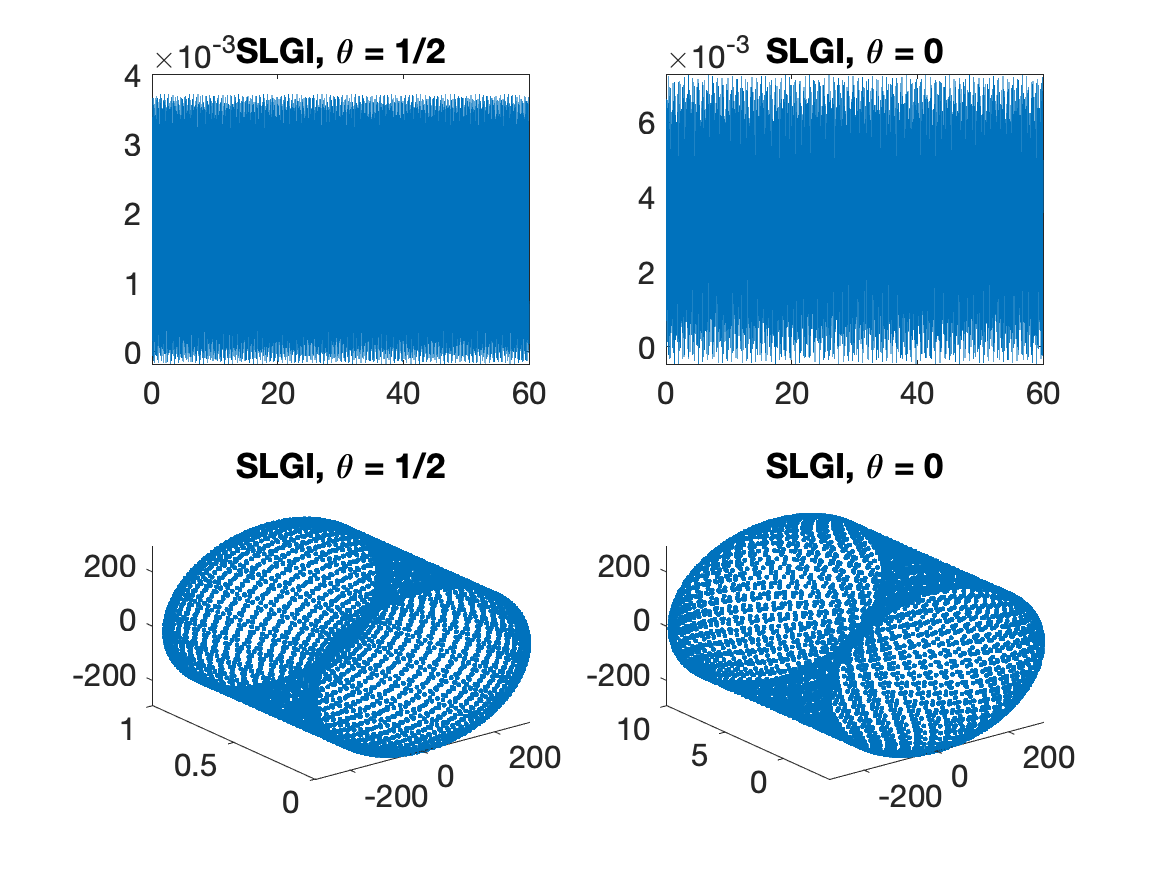}
    \caption{Symplectic Lie group integrators, long time integration, $h=0.01$, $6000$ steps.. Top: energy error, bottom 3D plot of  $M\ell Q^{-1}\Gamma_0$.}
    \label{fig:HT2}
\end{figure}

In Figure~\ref{fig:HT1} we report the performance of the \RE{symplectic} Lie group integrators \eqref{symplectic0}-\eqref{symplectic2} applied both on the equations \eqref{ht10}-\eqref{ht20} with $\theta=0$,  $\theta=\frac{1}{2}$ and $\theta=1$ (SLGI), and to the equations \eqref{ht211}-\eqref{ht241} with $\theta=\frac{1}{2}$ (SLGIKK). The methods with $\theta=\frac{1}{2}$ attain order $2$.
In Figure~\ref{fig:HT2} we show the energy error for the symplectic Lie group integrators with $\theta=\frac{1}{2}$ and $\theta=0$ integrating with stepsize $h=0.01$ for $6000$ steps. 

\section{Variable step size} \label{implementation}

One approach for varying the step size is based on the use of an embedded Runge--Kutta pair. This principle can be carried from standard Runge--Kutta methods in vector spaces to the present situation with RKMK and commutator-free schemes via minor modifications. We briefly summarise the main principle of embedded pairs before giving more specific details for the case of Lie group integrators. This approach is very well documented in the literature and goes back to Merson \cite{merson57aom} and a detailed treatment can be found in \cite[p. 165--168]{hairer93sod}.

An embedded pair consists of a main method used to propagate the numerical solution, together with some auxiliary method that is only used to obtain an estimate of the local error.
This local error estimate is in turn used to derive a step size adjustment formula that attempts to keep the local error estimate approximately equal to some user defined tolerance $\tol$ in every step.
Suppose the main method is of order $p$ and the auxiliary method is of order $\aux{p}\neq p.$\ \footnote{In this paper we will assume $\aux{p}<p$ in which case the local error estimate is relevant for the approximation $\aux{y}_{n+1}$} Both methods are applied to the input value $y_n$ and yields approximations $y_{n+1}$ and $\aux{y}_{n+1}$ respectively, using the same step size $h_{n+1}$. Now, some distance measure\footnote{\RE{There are many options for how to do this in practice, and the choice may also depend on the application. E.g. a Riemannian metric is a natural and robust alternative here.}} between 
$y_{n+1}$ and  $\aux{y}_{n+1}$ provides an estimate $e_{n+1}$ for the size of the local truncation error. Thus,
$e_{n+1}=Ch_{n+1}^{\aux{p}+1}+\mathcal{O}(h^{\aux{p}+2})$. Aiming at $e_{n+1}\approx\tol$ in every step, one may use a formula of the type
\begin{equation} \label{stepsizecontrol}
h_{n+1} = \theta\left(\frac{\tol}{\RE{e_{n+1}}}\right)^{\tfrac{1}{\aux{p}+1}}\, h_n
\end{equation}
where $\theta$ is a `safety factor', typically chosen between $0.8$ and $0.9$.
In case the step is rejected because $e_n>\tol$ we can redo the step with a step size obtained by the same formula. We summarise the approach in the following algorithm
\begin{tabbing}
xxx\=xxx\=xxx\=xxxx\kill\\
Given $y_n$, $h_n$, $\tol$ \\
Let $h:=h_n$\\
\textbf{repeat} \\
\>Compute $y_{n+1}$, $\aux{y}_{n+1}$, $e_{n+1}$ from $y_n$, $h$ \\
\>Update stepsize  $h:=\theta \left(\frac{\tol}{e_{n+1}}\right)^{\alpha} h$\\
\>accepted :=  $\RE{e_{n+1}}<\tol$ \\
\>\textbf{if} accepted\\
\>\>update step index: $n:=n+1$\\
\>\>$h_n:=h$ \\
\textbf{until} accepted
\end{tabbing}
Here we have used again the safety factor $\theta$, and the parameter $\alpha$ is generally chosen as $\alpha=\frac1{1+\min(p,\aux{p})}$.

\subsection{RKMK methods with variable stepsize}
We need to specify how to calculate the quantity $e_{n+1}$ in each step. For RKMK \RE{methods} the situation is simplified by the fact that we are solving the local problem \eqref{dexpinveq} in the linear space $\g$, where the known theory can be applied directly.
So any standard embedded pair of Runge--Kutta methods described by coefficients $(a_{ij}, b_i, \aux{a}_{ij}, \aux{b}_i)$ of orders $(p, \aux{p})$ can be applied to the full dexpinv-equation \eqref{dexpinveq} to obtain local Lie algebra approximations $\sigma_1$, $\aux{\sigma}_1$ and \RE{one uses} e.g.
$e_{n+1}=\|\sigma_1-\aux{\sigma}_1\|$ (note that the equation itself depends on $y_n$).
For methods which use a truncated version of the series for $\dexp_u^{-1}$ one may also try to optimise performance by including commutators that are shared between the main method and the auxiliary scheme.

\subsection{Commutator-free methods with variable stepsize}
For the commutator-free methods of section~\ref{twoclasses} the situation is different since such methods do not have a natural local representation in a linear space. One can still derive embedded pairs, and this can be achieved by studying order conditions \cite{owren06ocf} as was done in \cite{curry19vss}. Now one obtains after each step two approximations $y_{n+1}$ and $\aux{y}_{n+1}$ on $\M$ both by using the same initial value $y_n$ and step size $\RE{h_n}$. One must also have access to some metric $d$ to calculate $e_{n+1}=d(y_{n+1},\aux{y}_{n+1})$
We give a few examples of embedded pairs.

\subsubsection{Pairs of order $(p,\aux{p})=(3,2)$}

It is possible to obtain embedded pairs of order 3(2) which satisfy the  requirements above.
We present two examples from \cite{curry19vss}. The first one reuses the second stage exponential in the update
\addtolength{\jot}{1mm}
\begin{align*}
Y_{n,1} &= y_n \\
Y_{n,2} &= \exp(\tfrac13 hf_{n,1})\cdot y_n \\
Y_{n,3} &= \exp(\tfrac23 h f_{n,2})\cdot y_n \\
y_{n+1} &= \exp(h(-\tfrac{1}{12}f_{n,1}+\tfrac34 f_{n,3}))\cdot Y_{n,2} \\
\aux{y}_{n+1} &= \exp(\tfrac{1}{2} h (f_{n,2}+f_{n,3}))\cdot y_n
\end{align*}
One could also have reused the third stage $Y_{n,3}$ in the update, rather than $Y_{n,2}$. 
\begin{align*}
Y_{n,1} &= y_n \\
Y_{n,2} &= \exp(\tfrac23 hf_{n,1})\cdot y_n \\
Y_{n,3} &= \exp(h (\tfrac{5}{12}f_{n,1} + \tfrac14 f_{n,2})\cdot y_n \\
y_{n+1} &= \exp(h(-\tfrac{1}{6}f_{n,1}-\tfrac12 f_{n,2}+f_{n,3}))\cdot Y_{n,3} \\
\aux{y}_{n+1} &= \exp(\tfrac{1}{4} h (f_{n,1}+3f_{n,3}))\cdot y_n
\end{align*}
It is always understood that the frozen vector fields are $f_{n,i}:= f_{Y_{n,i}}$.

\subsubsection{Order $(4,3)$}\label{ssec:43}
The procedure of deriving efficient pairs \RE{becomes} more complicated as the order increases. In \cite{curry19vss} \RE{a low cost} pair of order $(4,3)$ was derived, in the sense that one attempted to minimise the number of stages \RE{and exponentials} in the embedded pair as a whole. This came, however, at the expense of a relatively large error constant. So rather than presenting the method from that paper, we suggest a simpler procedure at the cost of some more computational work per step, we simply furnish the commutator-free method of section~\ref{LGI} by a third order auxiliary scheme.
It can be described as follows:
\begin{enumerate}
    \item Compute $Y_{n,i},\;i=1\ldots,4$ and $y_{n+1}$ from \eqref{cf4} 
    \item Compute an additional stage $\bar{Y}_{n,3}$ and then $\aux{y}_{n+1}$ as
    \begin{align} \label{cf4e}
\begin{split}
\bar{Y}_{n,3} &= \exp(\tfrac34h f_{n,2})\cdot y_n \\
\aux{y}_{n+1} &= \exp(\tfrac{h}{9}(-f_{n,1}+3f_{n,2}+4\bar{f}_{n,3}))\cdot \exp(\tfrac{h}{3}f_{n,1}) \cdot y_n
\end{split}
\end{align}
\end{enumerate}

\section{The $N$-fold 3D pendulum}~\label{pendula}
In this section, we present a model for a system of $N$ connected 3-dimensional pendulums. The modelling part comes from \cite{lee18gfo}, and here we study the vector field describing the dynamics, in order to re-frame it into the Lie group integrators setting described in the previous sections.
The model we use is not completely realistic since, for example, it neglects possible interactions between pendulums, and it assumes ideal spherical joints between them. However, this is still a relevant example from the point of view of geometric numerical integration. More precisely, we show a possible way to work with a configuration manifold which is not a Lie group, applying the theoretical instruments introduced before. \begin{figure}[htbp]
    \centering
    \includegraphics[width=.6\textwidth]{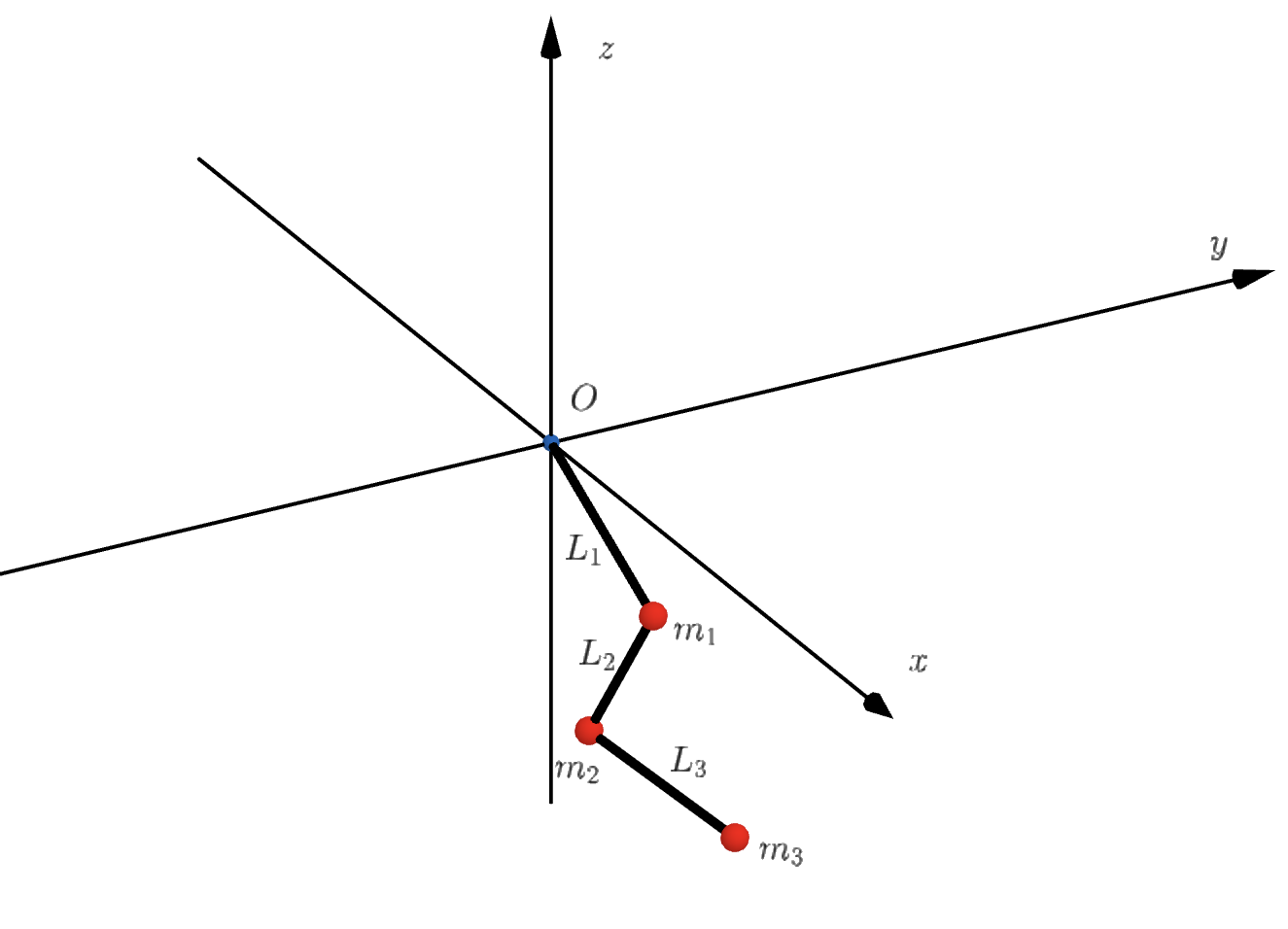}
    \caption{\footnotesize{$3-$fold pendulum at a fixed time instant, with fixed point placed at the origin.}}
    \label{fig:nfold}
\end{figure} 
The Lagrangian we consider is a function from $(TS^2)^N$ to $\mathbb{R}$. Instead of the coordinates $(q_1,...,q_N,\dot{q}_1,...,\dot{q}_N)$, where $\dot{q}_i\in T_{q_i}S^2$, we choose to work with the angular velocities. Precisely, 
$$
T_{q_i}S^2 = \{v\in\mathbb{R}^3:\;v^Tq_i=0\} = \langle q_i\rangle ^{\perp} \subset \mathbb{R}^3,
$$
and hence for any $\dot{q}_i\in T_{q_i}S^2$ there exist $\omega_i\in\mathbb{R}^3$ such that $\dot{q}_i=\omega_i\times q_i$, which can be interpreted as the angular velocity of $q_i$. So we can assume without loss of generality that $\omega_i^Tq_i=0$ (i.e. $\omega_i\in T_{q_i}S^2$) and pass to the coordinates $(q_1,\omega_1,q_2,\omega_2,...,q_N,\omega_N)\in (TS^2)^N$ to describe the dynamics.  In this section we denote with $m_1,...,m_N$ the masses of the pendulums and with $L_1,..., L_N$ their lengths. Figure \ref{fig:nfold} shows the case $N=3$.
We organize the section into three parts:
\begin{enumerate}
    \item We define the transitive Lie group action used to integrate this model numerically,
    \item We show a possible way to express the dynamics in terms of the infinitesimal generator of this action, for the general case of $N$ joint pendulums,
    \item We focus on the case $N=2$, as a particular example. For this setting, we present some numerical experiment comparing various Lie group integrators and some classical numerical integrator. Then we conclude with numerical experiments on variable step size.
\end{enumerate}
\subsection{Transitive group action on $(TS^2)^N$}\label{se:uno}

We characterize a transitive action for $(TS^2)^N$, starting with the case $N=1$ and generalizing it to $N>1$ . The action we consider is based on the identification between $\mathfrak{se}(3)$, the Lie algebra of $SE(3)$, and $\mathbb{R}^6$. We start from the Ad-action of $SE(3)$ on $\mathfrak{se}(3)$ (see \cite{Holm}), which writes
\[
\Ad : SE(3)\times \mathfrak{se}(3) \rightarrow \mathfrak{se}(3),
\]
\[
\Ad((R,r),(u,v)) = (Ru,Rv+\hat{r}Ru).
\]
Since $\mathfrak{se}(3)\simeq \mathbb{R}^6$, the Ad-action allows us to define the following Lie group action on $\mathbb{R}^6$
\[
\psi: SE(3)\times\mathbb{R}^6\rightarrow \mathbb{R}^6,\;\;\psi((R,r),(u,v)) = (Ru,Rv+\hat{r}Ru).
\]
We can think of $\psi$ as a Lie group action on $TS^2$ since, for any $q\in\mathbb{R}^3$, it maps points of
\[TS_{|q|}^2:=\{(\Tilde{q},\Tilde{\omega})\in \mathbb{R}^3\times\mathbb{R}^3:\; \Tilde{\omega}^T\Tilde{q}=0,\;|\Tilde{q}|=|q|\}\subset \mathbb{R}^6
\]
into other points of $TS_{|q|}^2$. Moreover, with standard arguments (see \cite{olver2000applications}), it is possible to prove that the orbit of a generic point $m=(q,\omega)\in\mathbb{R}^6$ with $\omega^Tq=0$ coincides with
\[
\RE{\text{Orb}}(m)=TS_{|q|}^2.
\]
In particular, when $q\in\mathbb{R}^3$ is a unit vector (i.e. $q\in S^2$), $\psi$ allows us to define a transitive Lie group action on $TS^2=TS_{|q|=1}^2$ which writes
\[
\psi : SE(3)\times TS^2 \rightarrow TS^2
\]
\[
\psi((A,a),(q,\omega)) := \psi_{(A,a)}(q,\omega) =  (Aq,A\omega + \hat{a}Aq)=(\bar{q},\bar{\omega}).
\]
To conclude the description of the action, we report here its infinitesimal generator which is fundamental in the Lie group integrators setting
\[
\left.\infgen((u,v))\right|_{(q,\omega)} =(\hat{u}q,\hat{u}\omega + \hat{v}q).
\]
We can extend this construction to the case $N>1$ in a natural way, i.e. through the action of a Lie group obtained from cartesian products of $SE(3)$ and equipped with the direct product structure. More precisely, we consider the group $G=(SE(3))^N$ and by direct product structure we mean that for any pair of elements $$\delta^{(1)}=(\delta^{(1)}_1,...,\delta^{(1)}_N),\quad \delta^{(2)}=(\delta^{(2)}_1,...,\delta^{(2)}_N)\in G,$$ denoted with $*$ the semidirect product of $SE(3)$, we define the product $\circ$ on $G$ as
\[
\delta^{(1)}\circ \delta^{(2)} := (\delta^{(1)}_1 * \delta^{(2)}_1,...,\delta^{(1)}_N * \delta^{(2)}_N)\in G.
\]
With this group structure defined, we can generalize the action introduced for $N=1$ to larger $N$s as follows
\[
\psi : (SE(3))^N\times (TS^2)^N \rightarrow (TS^2)^N,
\]
\[
\begin{split}
\psi&((A_1,a_1,...,A_N,a_n),(q_1,\omega_1,...,q_N,\omega_N)) =\\ &=(A_1q_1,A_1\omega_1+\hat{a}_1A_1q_1,...,A_Nq_N,A_N\omega_N+\hat{a}_NA_Nq_N),
\end{split}
\]
whose infinitesimal generator writes
\[
\infgen(\xi)\vert_m =(\hat{u}_1q_1,\hat{u}_1\omega_1+\hat{v}_1q_1,...,\hat{u}_Nq_N,\hat{u}_N\omega_N+\hat{v}_Nq_N),
\]
where $\xi=[u_1,v_1,...,u_N,v_N]\in\mathfrak{se}(3)^N$ and $m=(q_1,\omega_1,...,q_N,\omega_N)\in (TS^2)^N$.
We have now the only group action we need to deal with the $N-$fold spherical pendulum. In the following part of this section we work on the vector field describing the dynamics and adapt it to the Lie group integrators setting.

\subsection{Full chain}
We consider the vector field $F\in\mathfrak{X}((TS^2)^N)$, describing the dynamics of the $N$-fold 3D pendulum, and we express it in terms of the infinitesimal generator of the action defined above. More precisely, we find a function $F:(TS^2)^N\rightarrow \mathfrak{se}(3)^N$ such that
\[
\infgen(f(m))\vert_m = F\vert_m,\;\;\forall m\in (TS^2)^N.
\]
We omit the derivation of $F$ starting from the Lagrangian of the system, which can be found in the section devoted to mechanical systems on $(S^2)^N$ of \cite{lee18gfo}. 
The configuration manifold of the system is $(S^2)^N$, while the Lagrangian, expressed in terms of the variables $(q_1,\omega_1,...,q_N,\omega_N)\in (TS^2)^N$, writes
\[
L(q,\omega) = T(q,\omega)-U(q) =\frac{1}{2}\sum_{i,j=1}^N\Big(M_{ij}\omega_i^T\hat{q}_i^T\hat{q}_j\omega_j\Big) - \sum_{i=1}^N\Big(\sum_{j=i}^N m_j\Big)gL_ie_3^Tq_i,
\]
where
\[
M_{ij} =\Big(\sum_{k=\RE{\text{max}}\{i,j\}}^N m_k\Big)L_iL_j I_3\in\mathbb{R}^{3\times 3}
\]
is the inertia matrix of the system\RE{, $I_3$ is the $3\times 3$ identity matrix,} and $e_3 = [0,0,1]^T$. Noticing that when $i=j$ we get
\[
\omega_i^T\hat{q}_i^T\hat{q}_i\omega_i = \omega_i^T(I_3-q_iq_i^T)\omega_i = \omega_i^T\omega_i,
\]
we simplify the notation writing 
\[
T(q,\omega) = \frac{1}{2}\sum_{i,j=1}^N\Big(\omega_i^TR(q)_{ij}\omega_j\Big)
\]
where $R(q)\in\mathbb{R}^{3N\times 3N}$ is a symmetric block matrix defined as
\[
R(q)_{ii} = \Big(\sum_{j=i}^Nm_j\Big)L_i^2I_3\in\mathbb{R}^{3\times 3},
\]
\[
R(q)_{ij} = \Big(\sum_{k=j}^N m_k\Big)L_iL_j\hat{q}_i^T\hat{q}_j\in\mathbb{R}^{3\times 3} = R(q)_{ji}^T,\; i<j.
\]
The vector field on which we need to work defines the following first-order ODE
\[
\begin{split}
\dot{q}_i &= \omega_i \times q_i,\;i=1,...,N,\\
R(q)\dot{\omega} &=\left[\sum_{\substack{j=1\\ j\neq i}}^N M_{ij}|\omega_j|^2\hat{q}_iq_j - \Big(\sum_{j=i}^N m_j\Big)gL_i\hat{q}_ie_3 \right]_{i=1,...,N}\in\mathbb{R}^{3N}
\end{split}
\]
By direct computation it is possible to see that, for any $q=(q_1,...,q_N)\in (S^2)^N$ and $\omega\in T_{q_1}S^2\times ... \times T_{q_N}S^2$, we have
\[
(R(q)\omega)_i \in T_{q_i}S^2.
\]
Therefore, there is a well-defined linear map
\[
A_{q}:T_{q_1}S^2\times ... \times T_{q_N}S^2 \rightarrow T_{q_1}S^2 \times ... \times T_{q_N}S^2, A_q(\omega):=R(q)\omega.
\]
We can even notice that $R(q)$ defines a positive-definite bilinear form on this linear space, since
\[
\omega^TR(q)\omega = \sum_{i,j=1}^N \omega_i^T\hat{q}_i^TM_{ij}\hat{q}_j\omega_j = \sum_{i,j=1}^N (\hat{q}_i\omega_i)^TM_{ij}(\hat{q}_j\omega_j) = v^TMv>0.
\]
The last inequality holds because $M$ is the inertia matrix of the system and hence it defines a symmetric positive-definite bilinear form on $T_{q_1}S^2\times...\times T_{q_N}S^2$, see e.g. \cite{Goldstein}  \footnote{It follows from the definition of the inertia tensor, i.e. $$0\leq \Tilde{T}(q,\dot{q}) =\frac{1}{2} \sum_{i=1}^N \Big(\sum_{j\geq i}m_j\Big)L_iL_j\dot{q}_i^T\dot{q}_j := \frac{1}{2}\dot{q}^TM\dot{q}.$$ Moreover, in this situation it is even possible to explicitly find the Cholesky \RE{factorization of the matrix} $M$ with an iterative algorithm.}. This implies the map $A_q$ is invertible and hence we are ready to express the vector field in terms of the infinitesimal generator. We can rewrite the ODEs for the angular velocities as follows:
\[
\dot{\omega}= A_{q}^{-1}\Big([ g_1 , ... , g_N  ]^T\Big) =\begin{bmatrix}
h_1(q,\omega) \\ ... \\ h_N(q,\omega)
\end{bmatrix} = \begin{bmatrix}
a_1(q,\omega)\times q_1 \\
...\\
a_N(q,\omega)\times q_N
\end{bmatrix}
\]
where
\[
g_i=g_i(q,\omega) = \sum_{\substack{j=1\\ j\neq i}}^N M(q)_{ij}|\omega_j|^2\hat{q}_iq_j - \Big(\sum_{j=i}^N m_j\Big)gL_i\hat{q}_ie_3,\;i=1,...,N
\]
and $a_1,...,a_N:(TS^2)^N\rightarrow \mathbb{R}^3$ are $N$ functions whose existence is guaranteed by the analysis done above. Indeed, we can set $a_i(q,\omega):=q_i\times h_i(q,\omega)$ and conclude that a mapping $f$ from $(TS^2)^N$ to $(\mathfrak{se}(3))^N$ such that
\[
\infgen(f(q,\omega))\vert_{(q,\omega)}=F\vert_{(q,\omega)}
\]
is the following one
\[
f(q,\omega) = \begin{bmatrix}
\omega_1 \\
q_1\times h_1 \\
...\\
...\\
\omega_N \\
q_N\times h_N
\end{bmatrix}\in\mathfrak{se}(3)^N\simeq \mathbb{R}^{6N}.
\]
We will not go into the Hamiltonian formulation of this problem; however, we remark that a similar approach works even in that situation. Indeed, following the derivation presented in \cite{lee18gfo}, we see that for a mechanical system on $(S^2)^N$ the conjugate momentum writes
\[
T_{q_1}^*S^2\times ... T_{q_N}^*S^2\ni\pi=(\pi_1,...,\pi_N),\text{ where }\pi_i = -\hat{q}_i^2\frac{\partial L}{\partial \omega_i}
\]
and its components are still orthogonal to the respective base points $q_i\in S^2$. Moreover, Hamilton's equations take the form
\[
\begin{split}
\dot{q}_i &= \frac{\partial H(q,\pi)}{\partial \pi_i}\times q_i, \\
\dot{\pi}_i &= \frac{\partial H(q,\pi)}{\partial q_i}\times q_i + \frac{\partial H(q,\pi)}{\partial \pi_i}\times \pi_i,
\end{split}
\]
which implies that setting 
\[ f(q,\pi) = \begin{bmatrix} \partial_{q_1}H(q,\pi), & \partial_{\pi_1}H(q,\pi), & \dots ,& \partial_{q_N}H(q,\pi), & \partial_{\pi_N}H(q,\pi)\end{bmatrix}
\]
we can represent even the Hamiltonian vector field of the $N-$fold 3D pendulum in terms of this group action.

\subsubsection{Case $N=2$}

We have seen how it is possible to turn the equations of motion of a $N-$chain of pendulums into the Lie group integrators setting. Now we focus on the example with $N=2$ pendulums. The equations of motion write
\[
\dot{q}_1 = \hat{\omega}_1q_1,\quad \dot{q}_2 = \hat{\omega}_2q_2,
\]
\begin{equation}\label{eq:double}
R(q)\begin{bmatrix}
\dot{\omega}_1 \\ \dot{\omega}_2
\end{bmatrix}= 
\begin{bmatrix}
(-m_2L_1L_2|\omega_2|^2\hat{q}_2 + (m_1+m_2)gL_1\hat{e}_3)q_1 \\
(-m_2L_1L_2|\omega_1|^2\hat{q}_1 + m_2gL_2\hat{e}_3)q_2
\end{bmatrix},
\end{equation}
where 
\[
R(q) = \begin{bmatrix}
(m_1+m_2)L_1^2I_3 & m_2L_1L_2\hat{q}_1^T\hat{q}_2 \\
m_2L_1L_2\hat{q}_2^T\hat{q}_1 & m_2L_2^2I_3
\end{bmatrix}.
\]
As presented above, the matrix $R(q)$ defines a linear invertible map of the space $T_{q_1}S^2\times T_{q_2}S^2$ onto itself:
\[
A_{(q_1,q_2)}:T_{q_1}S^2\times T_{q_2}S^2\rightarrow T_{q_1}S^2\times T_{q_2}S^2,\;[\omega_1,\omega_2]^T\rightarrow R(q)[\omega_1,\omega_2]^T.
\]
We can easily see that it is well defined since
\[
R(q)\begin{bmatrix}
\omega_1 \\ \omega_2
\end{bmatrix} = \begin{bmatrix}
(m_1+m_2)L_1^2I_3 & m_2L_1L_2\hat{q}_1^T\hat{q}_2 \\
m_2L_1L_2\hat{q}_2^T\hat{q}_1 & m_2L_2^2I_3
\end{bmatrix}\begin{bmatrix}
\hat{v}_1q_1 \\ \hat{v}_2q_2
\end{bmatrix} = \begin{bmatrix}
\hat{r}_1q_1\\ \hat{r}_2q_2 
\end{bmatrix}\in (TS^2)^2
\]
with $$r_1(q,\omega):=(m_1+m_2)L_1^2v_1+m_2L_1L_2\hat{q}_2\hat{v}_2q_2,$$ $$ r_2(q,\omega):=m_2L_1L_2\hat{q}_1\hat{v}_1q_1+m_2L_2^2v_2. $$
This map guarantees that if we rewrite the pair of equations for the angular velocities in (\ref{eq:double}) as
\[
\begin{split}
\dot{\omega}&= R^{-1}(q)\begin{bmatrix}
(-m_2L_1L_2|\omega_2|^2\hat{q}_2 + (m_1+m_2)gL_1\hat{e}_3)q_1 \\
(-m_2L_1L_2|\omega_1|^2\hat{q}_1 + m_2gL_2\hat{e}_3)q_2
\end{bmatrix}=R^{-1}(q)b=\\
&=A_{(q_1,q_2)}^{-1}(b)=\begin{bmatrix}
h_1 \\ h_2
\end{bmatrix}\in T_{q_1}S^2\times T_{q_2}S^2,
\end{split}
\]
then we are assured \RE{that there exists a pair of functions} $a_1,a_2:TS^2\times TS^2\rightarrow\mathbb{R}^3$ such that
\[
\dot{\omega} = \begin{bmatrix}
a_1(q,\omega)\times q_1 \\ a_2(q,\omega)\times q_2
\end{bmatrix} = \begin{bmatrix}
h_1(q) \\ h_2(q)
\end{bmatrix}.
\]
Since we want $a_i\times q_i = h_i$, we just impose $a_i=q_i\times h_i$ and hence the whole vector field can be rewritten as
\[
\begin{bmatrix}
\dot{q}_1 \\ \dot{\omega}_1 \\ \dot{q}_2 \\ \dot{\omega}_2
\end{bmatrix} = \begin{bmatrix}
\omega_1 \times q_1 \\ (q_1\times h_1)\times q_1 \\ \omega_2\times q_2 \\ (q_2\times h_2)\times q_2
\end{bmatrix} = F\vert_{(q,\omega)},
\]
with $h_i=h_i(q,\omega)$ and
\[
\begin{bmatrix}
h_1(q,\omega) \\ h_2(q,\omega)
\end{bmatrix} = R^{-1}(q)\begin{bmatrix}
(-m_2L_1L_2|\omega_2|^2\hat{q}_2 + (m_1+m_2)gL_1\hat{e}_3)q_1 \\
(-m_2L_1L_2|\omega_1|^2\hat{q}_1 + m_2gL_2\hat{e}_3)q_2
\end{bmatrix}.
\]
Therefore, we can express the whole vector field in terms of the infinitesimal generator of the action of $SE(3)\times SE(3)$ as
\[
\infgen(f(q,\omega))\vert_{(q,\omega)}=F\vert_{(q,\omega)}
\]
through the function
\[
f : TS^2\times TS^2\rightarrow \mathfrak{se}(3)\times\mathfrak{se}(3)\simeq \mathbb{R}^{12},\;\;(q,\omega)\rightarrow (\omega_1, q_1\times h_1, \omega_2,q_2\times h_2).\]

\subsection{Numerical experiments}
In this section, we present some numerical experiment for the $N-$chain of pendulums. We start by comparing the various Lie group integrators that we have tested (with the choice $N=2$), and conclude by analyzing an implementation of variable step size. Lie group integrators allow to keep the evolution of the solution in the correct manifold, which is $TS^2\times TS^2$ when $N=2$. Hence, we briefly report two sets of numerical experiments. In the first one, we show the convergence rate of all the Lie group integrators tested on this model. In the second one, we check how they behave in terms of preserving the two following relations:
\begin{itemize}
     \item $q_i(t)^T q_i(t) = 1,\text{ i.e. } q_i(t)\in S^2,\;i=1,2,$
    \item $q_i(t)^T\omega_i(t) = 0,\text{ i.e. }\omega_i(t)\in T_{q_i(t)}S^2,\;i=1,2,$
\end{itemize}
completing the analysis with a comparison with the classical Runge--Kutta 4 and with ODE45 of MATLAB. The Lie group integrators used to obtain the following experiments are Lie Euler, Lie Euler Heun, three versions of  Runge--Kutta--Munthe--Kaas methods of order four and one of order three. The RKMK4 with two commutators mentioned in the plots, is the one presented in Section \ref{LGI}, while the other schemes can be found for example in \cite{celledoni14ait}. 

\begin{figure}[htbp]
    \centering
    \includegraphics[trim= 400 0 400 0,width=.8\textwidth]{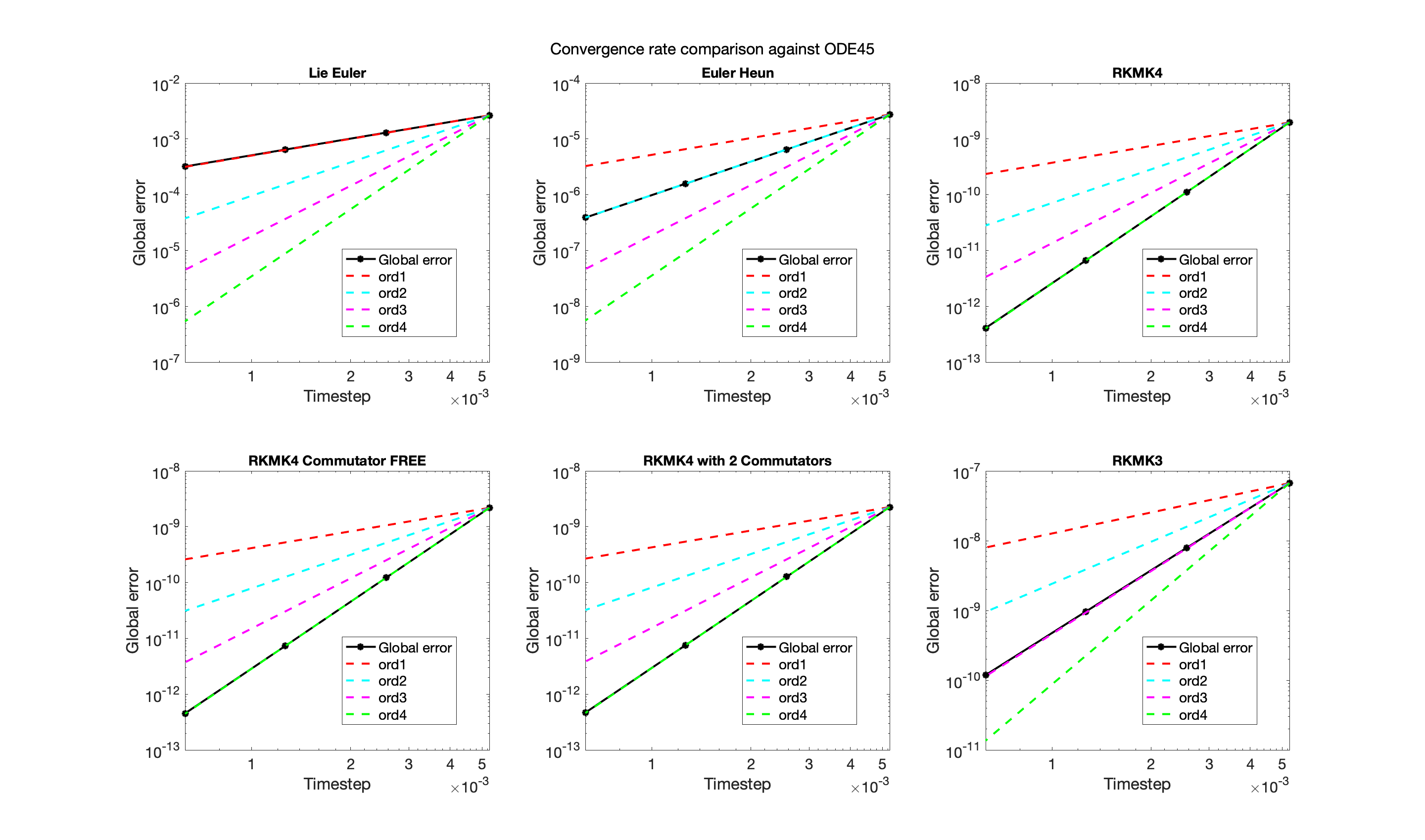}
    \caption{Convergence rate of the implemented Lie group integrators, based on global error considering as a reference solution the one of ODE45, with strict tolerance.}
    \label{fig:norms}
\end{figure}

Figure \ref{fig:norms} presents the plots of the errors, in logarithmic scale, obtained considering as a reference solution the one given by the ODE45 method, with strict tolerance. Here, we used an exact expression for the $\dexp_{\sigma}^{-1}$ function. However, we could obtain the same results with a truncated version of this function, keeping a sufficiently high number of commutators\RE{, or after some clever manipulations of the commutators (as with RKMK4 with 2 commutators, see Section \ref{twoclasses})}. The schemes show the right convergence rates, so we can move to the analysis of the time evolution on $TS^2\times TS^2$.

\begin{figure}[htbp]
    \centering
    \includegraphics[trim=200 0 200 0,width=\textwidth]{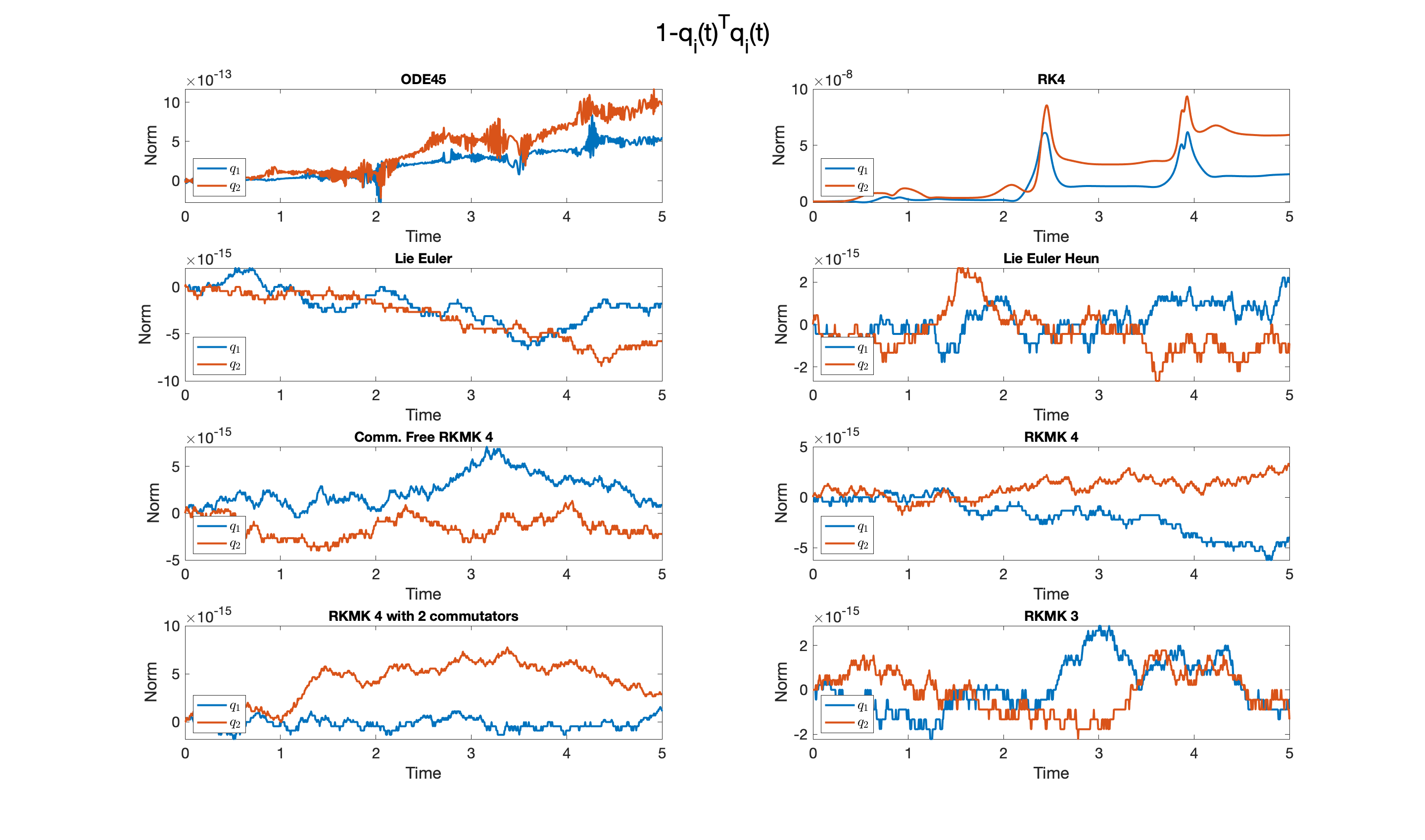}
    \caption{Visualization of the quantity $1-q_i(t)^Tq_i(t)$, $i=1,2$, for time $t\in [0,5]$. These plots focus on the preservation of the geometry of $S^2$.}
    \label{fig:S2}
\end{figure}

\begin{figure}[htbp]
    \centering
    \includegraphics[trim= 200 0 200 0,clip,width=\textwidth]{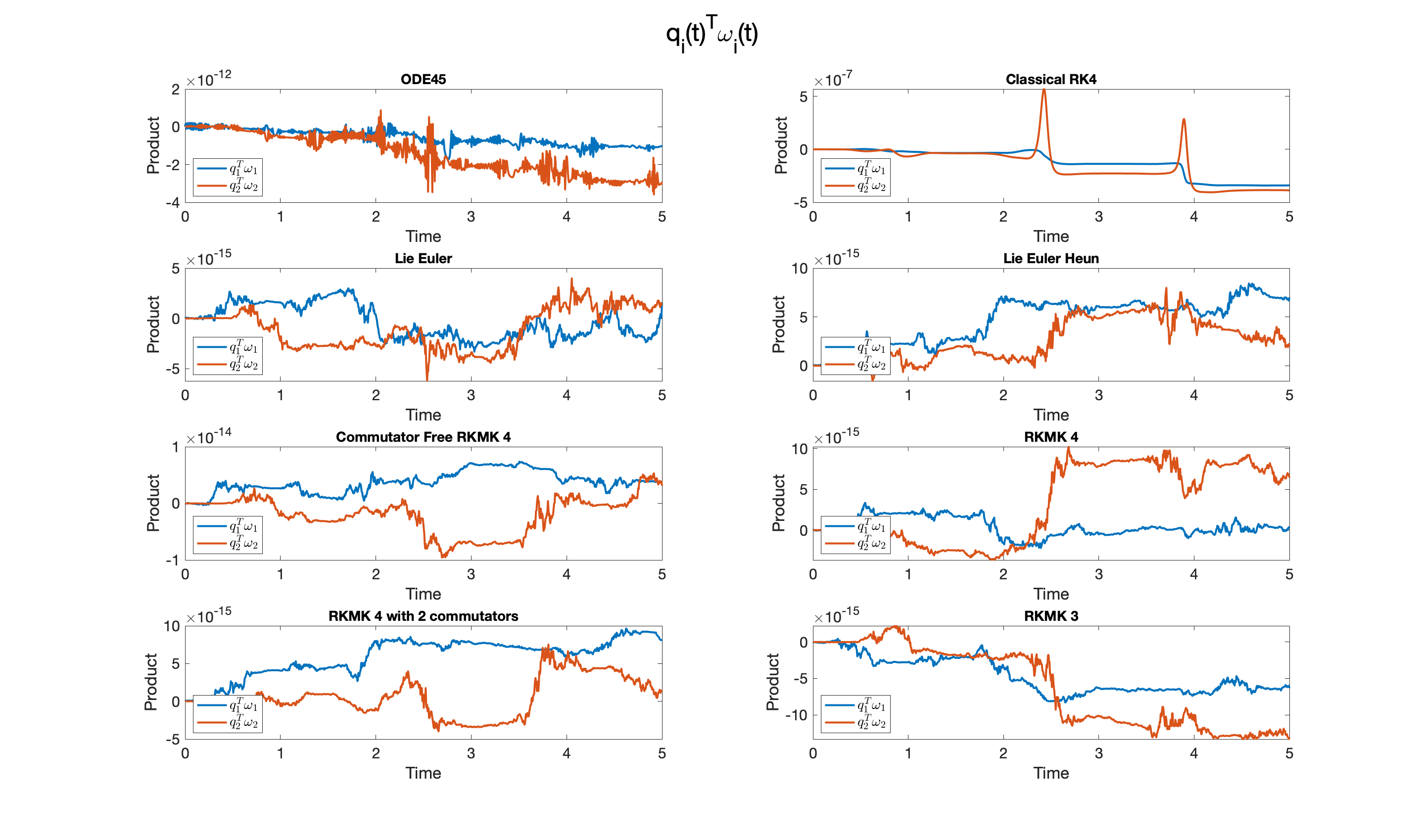}
    \caption{Visualization of the inner product $q_i(t)^T \omega_i(t)$, $i=1,2$, for $t\in [0,5]$. These plots focus on the preservation of the geometry of $T_{q_i(t)}S^2$.}
    \label{fig:tangent}
\end{figure}

In Figure \ref{fig:S2} we can see the comparison of the time evolution of the $2-$norms of $q_1(t)$ and $q_2(t)$, for $0\leq t\leq T=5$. As highlighted above, unlike classical numerical integrators like the one implemented in ODE45 or the Runge--Kutta 4, the Lie group methods preserve the norm of the base components of the solutions, i.e. $|q_1(t)|=|q_2(t)|=1$ $\forall t\in [0, T]$. Therefore, as expected, these integrators preserve the configuration manifold. However, to complete this analysis, we show the plots making a similar comparison but with the tangentiality conditions.
Indeed, in Figure \ref{fig:tangent} we compare the time evolutions of the inner products $q_1(t)^T\omega_1(t)$ and $q_2(t)^T\omega_2(t)$ for $t\in [0,5]$, i.e. we see if these integrators preserve the geometry of the whole phase space $TS^2\times TS^2$. As we can see, while for Lie group methods these inner products are of the order of $10^{-14}$ and $10^{-15}$, the ones obtained with classical integrators show that the tangentiality conditions are not preserved with the same accuracy.

We now move to some experiments on variable stepsize. In this last part we focus on the RKMK pair coming from Dormand–Prince method (DOPRI 5(4) \cite{dormand1980family}), which we denote with RKMK(5,4). The aim of the plots we show is to compare the same schemes, both with constant and variable stepsize. We start by setting a tolerance and solving the system with the RKMK(5,4) scheme. Using the same number of time steps, we solve it again with RKMK of order 5. These experiments show that, for some tolerance and some initial conditions, the step size's adaptivity improves the numerical approximation accuracy. Since we do not have an available analytical solution to quantify these two schemes' accuracy, we compare them with the solution obtained with a strict tolerance and ODE45. \RE{We compute such accuracy, at time $T=3$, by means of the Euclidean norm of the ambient space $\mathbb{R}^{6N}$.}
\begin{figure}[htbp]
    \centering
    \includegraphics[trim= 300 0 300 0,width=.8\textwidth]{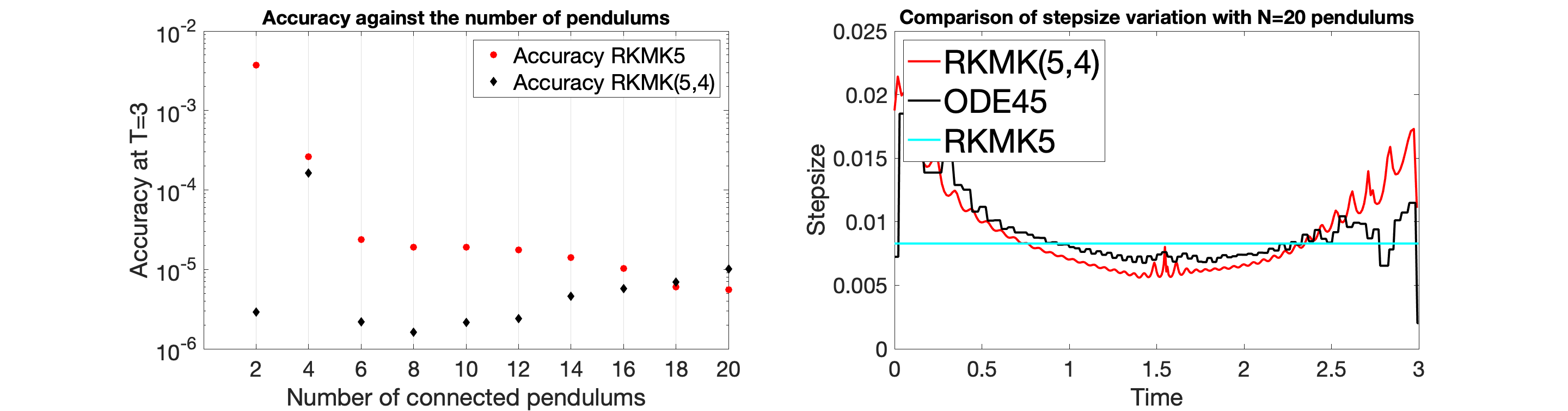}
    \caption{Comparison of accuracy at final time (on the left) and step adaptation for the case $N=20$ (on the right)\RE{, with all pendulums of length $L_i=1$}.}
    \label{fig:adapt}
\end{figure}

In Figure \ref{fig:adapt}, we compare the performance of the constant and variable stepsize methods, where the structure of the initial condition is always the same, but what changes is the number of connected pendulums. The considered initial condition is $ (q_i,\omega_i) = \left(\sqrt{2}/2,0,\sqrt{2}/2,0,1,0\right),\quad \forall i=1,...,N $, and all the masses and lengths are set to 1. From these experiments we can notice situations where the variable step size beats the constant one in terms of accuracy at the final time, like the case $N=2$ which we discuss in more detail afterwards. 

\smallskip

\RE{The results presented in Figure \ref{fig:adaptX} (left) do not aim to highlight any particular relation between how the number of pendulums increases or the regularity of the solution. Indeed, as we add more pendulums, we keep incrementing the total length of the chain since $\sum_{i=1}^NL_i=N$. Thus, here we do not have any appropriate limiting behaviour in the solution as $N\rightarrow +\infty$. The behaviour presented in that figure seems to highlight an improvement in accuracy for the RKMK5 method as $N$ increases. However, this is biased by the fact that when we increase $N$, to achieve the fixed tolerance of $10^{-6}$ with RKMKK(5,4), we need more time steps in the discretization. Thus, this plot does not say that as $N$ increases, the dynamics becomes more regular; it suggests that the number of required timesteps increases faster than the ``degree of complexity" of the dynamics.}\newline
\begin{figure}[htbp]
\centering
\begin{subfigure}{0.48\textwidth}
    \centering
    \includegraphics[trim=300 0 700 0,width=\textwidth]{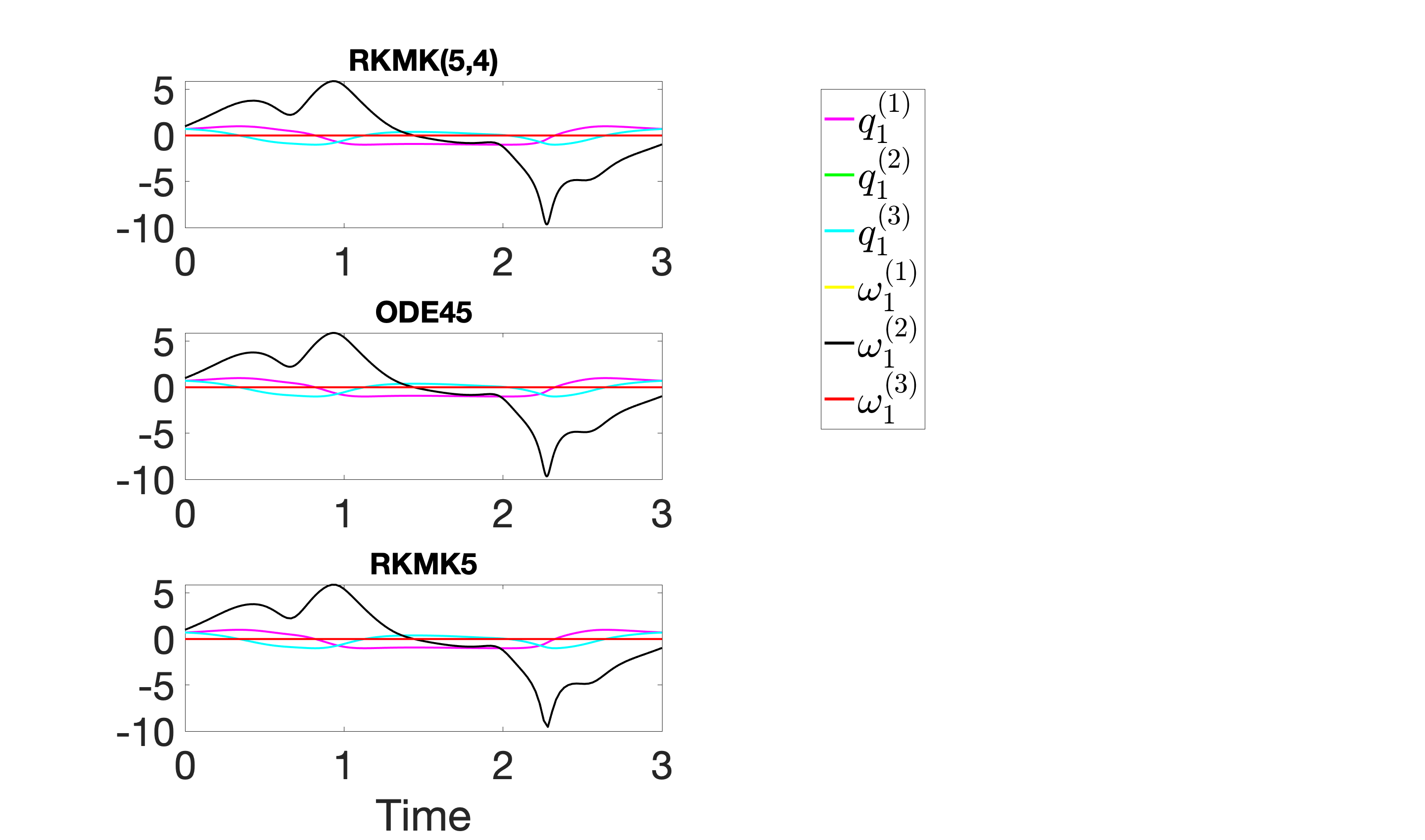}
    \caption{$(q_1(t),\omega_1(t))$}
    \label{fig:q1wi}
\end{subfigure}
\begin{subfigure}{0.48\textwidth}
    \centering
    \includegraphics[trim=300 0 700 0,width=\textwidth]{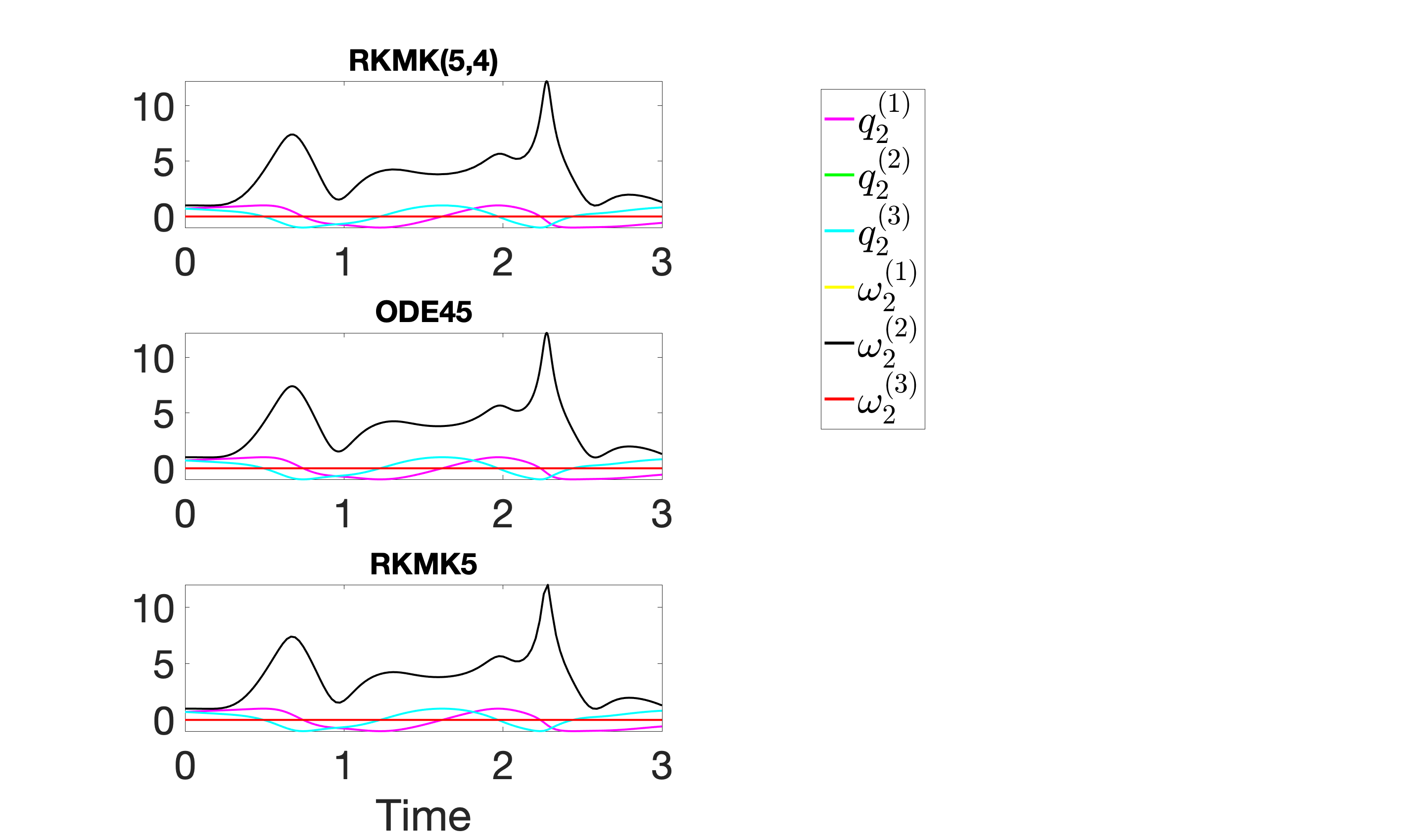}
    \caption{$(q_2(t),\omega_2(t))$}
    \label{fig:q2w2}
\end{subfigure}
\caption{In these plots we represent the six components of the solution describing the dynamics of the first mass (on the left) and of the second mass (on the right), for the case $N=2$. We compare the behaviour of the solution obtained with constant stepsize RKMK5, the variable stepsize RKMK(5,4) and ODE45.}
\label{fig:components}
\end{figure}

For the case $N=2$, we notice a relevant improvement passing to variable stepsize. In Figures \ref{fig:components} and \ref{fig:behaviour} we can see that, for this choice of the parameters, the solution behaves smoothly in most of the time interval, but then there is a peak in the second component of the angular velocities of both the masses, at $t\approx 2.2$. We can observe this behaviour both in the plots of Figure \ref{fig:components}, where we project the solution on the twelve components and even in Figure \ref{fig:variation}. In the latter, we plot two of the vector field components, i.e. the second components of the angular accelerations $\dot{\omega}_i(t)$, $i=1,2$. They show an abrupt change in the vector field in correspondence to $t\approx 2.2$, where the step is considerably restricted. \RE{Thus, to summarize, the gain we see with variable stepsize when $N=2$ is motivated by the unbalance in the length of the time intervals with no abrupt changes in the dynamics and those where they appear. Indeed, we see that apart from a neighbourhood of $t\approx 2.2$, the vector field does not change quickly. On the other hand, for the case $N=20$, this is the case. Thus, the adaptivity of the stepsize does not bring relevant improvements in the latter situation. } 

\smallskip

\RE{The motivating application behind our choice of this mechanical system has been some intuitive relation with a beam model, as highlighted in the introduction of this work. However, for this limiting behaviour to make sense, we should fix the length of the entire chain of pendulums to some $L$ (the length of the beam at rest) and then set the size of each pendulum to $L_i=L/N$. In this case, keeping the same tolerance of $10^{-6}$ for RKMK(5,4), we get the results presented in the following plot. We do not investigate more in details this approach, which might be relevant for further work, however we highlight that here the step adaptivity improves the results as we expected.}
\begin{figure}[htbp]
    \centering
    \includegraphics[trim= 300 0 300 0,width=.8\textwidth]{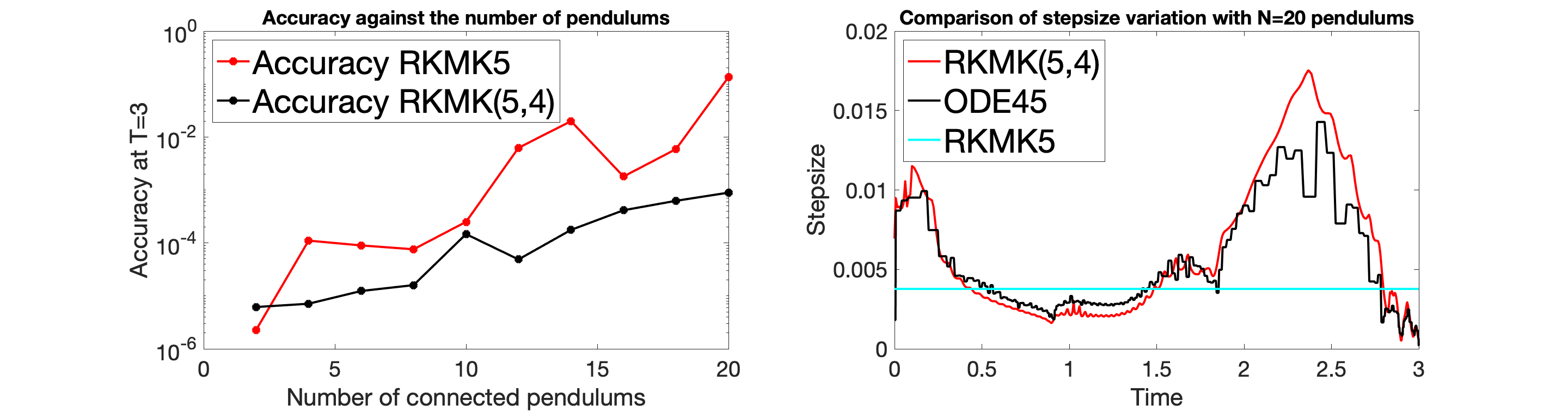}
    \caption{\RE{Comparison of accuracy at final time (on the left) and step adaptation for the case $N=20$ (on the right), with all pendulums of length $L_i=5/N$}.}
    \label{fig:adaptX}
\end{figure}

\begin{figure}[htbp]
\centering 
\begin{subfigure}{0.31\textwidth}
    \centering
    \includegraphics[trim=100 0 50 0,width=\textwidth]{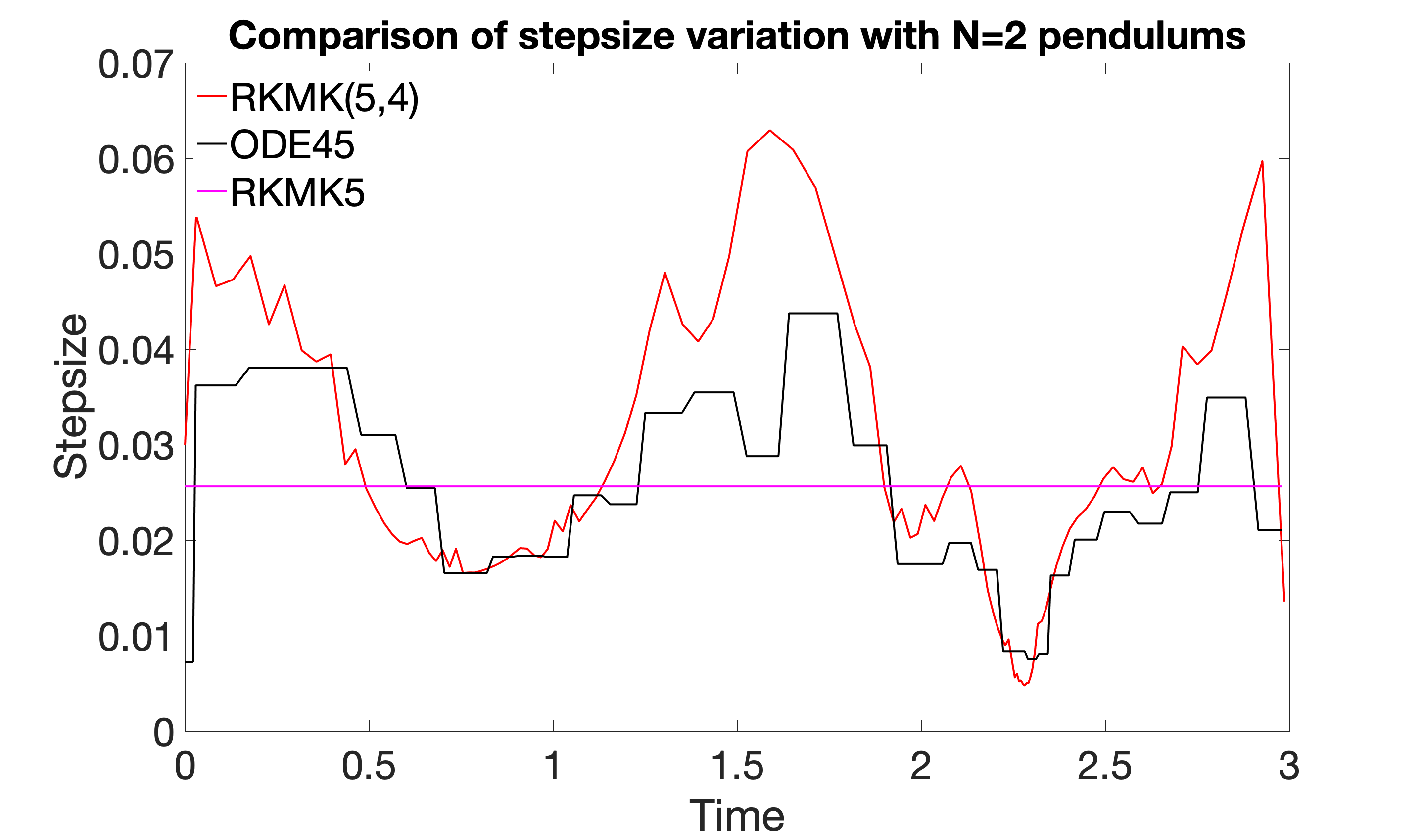}
    \caption{Step adaptation}
    \label{fig:adaptP2}
\end{subfigure}
\begin{subfigure}{0.33\textwidth}
    \centering
    \includegraphics[trim=110 0 50 0,width=\textwidth]{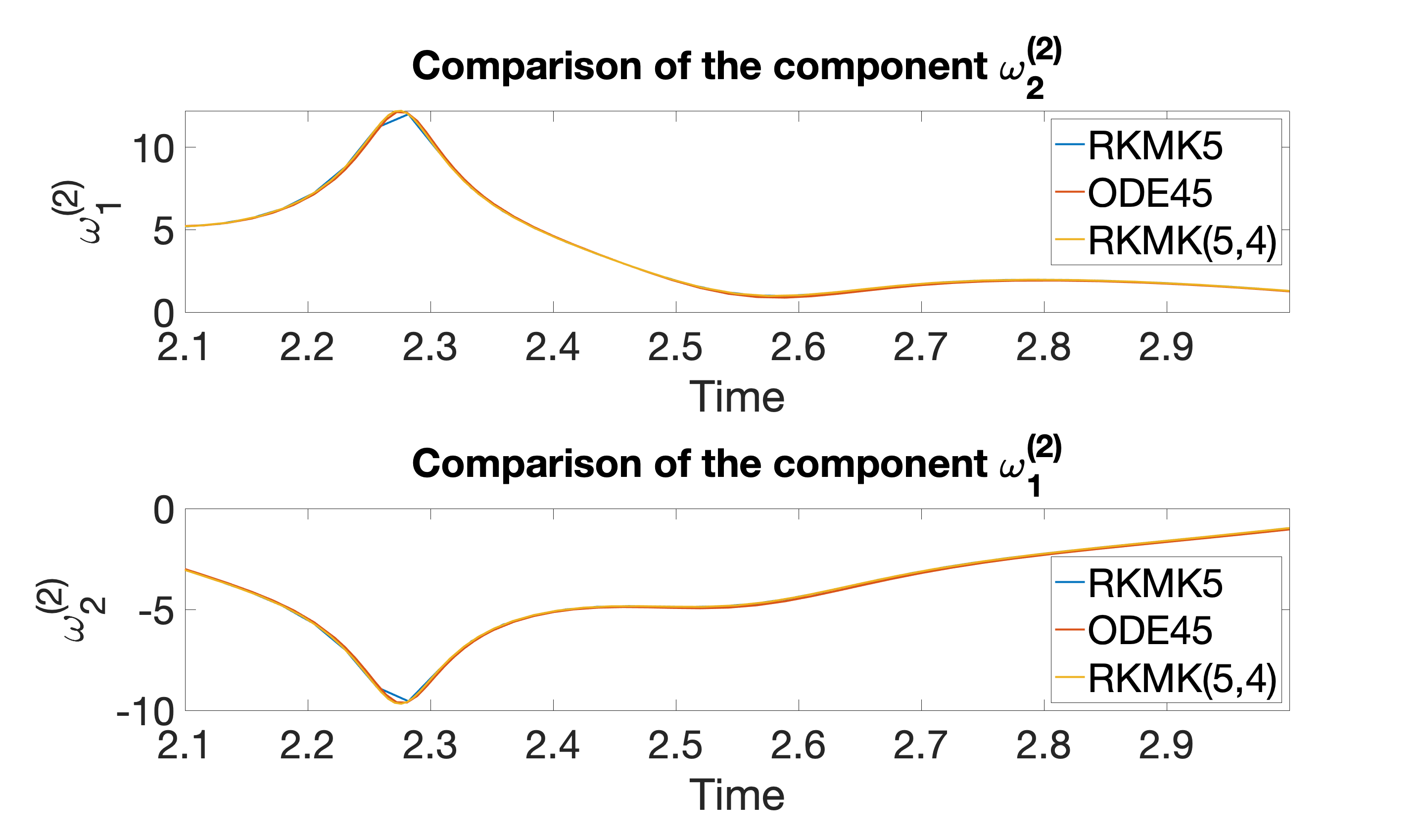}
    \caption{Zoom at final times}
    \label{fig:zoom}
\end{subfigure}
\begin{subfigure}{0.31\textwidth}
    \centering
    \includegraphics[trim=100 0 50 0,width=\textwidth]{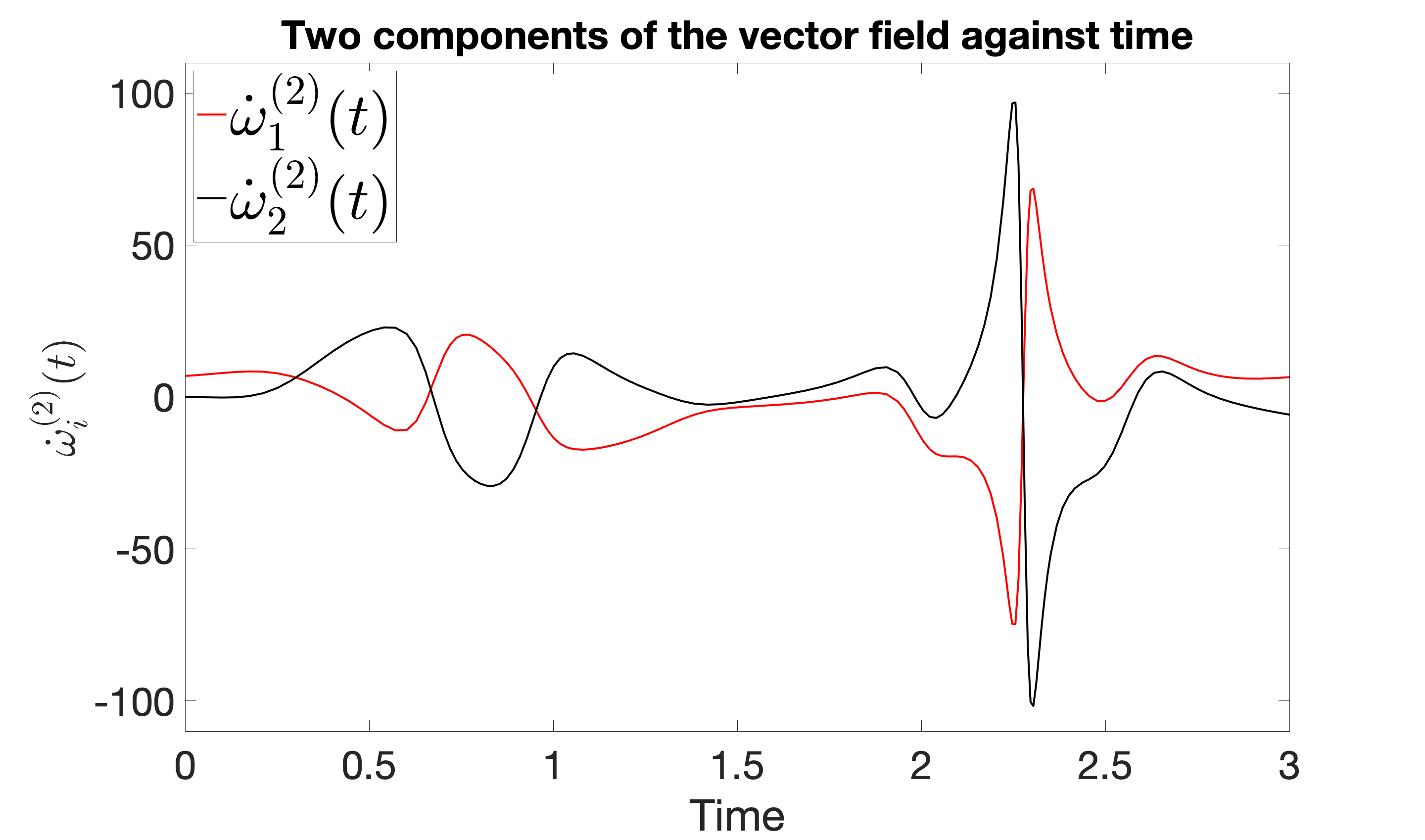}
    \caption{Values of $\dot{\omega}_i^{(2)}(t)$}
    \label{fig:variation}
\end{subfigure}
\caption{On the left, we compare the adaptation of the stepsize of RKMK(5,4) with the one of ODE45 and with the constant stepsize of RKMK5. In the center we plot the second component of the angular velocities $\omega_i^{(2)}$, $i=1,2$, and we zoom in the last time interval $t\in [2.1,3]$ to see that the variable stepsize version of the method better reproduces the reference solution. On the right, we visualize the speed of variation of second component of the angular velocities.}
\label{fig:behaviour}
\end{figure}

\section{Dynamics of two quadrotors transporting a mass point}~\label{quadrot}
In this section we consider a multibody system made of two cooperating quadrotor unmanned aerial vehicles (UAV) connected to a point mass (suspended load) via rigid links. This model is described in \cite{lee18gfo, Lee2013GeometricCO}.

We consider an inertial frame whose third axis goes in the direction of gravity, but opposite orientation, and we denote with $y\in\mathbb{R}^3$ the mass point and with $y_1,y_2\in\mathbb{R}^3$ the two quadrotors.
We assume that the links between the two quadrotors and the mass point are of a fixed length $L_1, L_2\in\mathbb{R}^+$. The configuration variables of the system are: the position of the mass point in the inertial frame, $y\in \mathbb{R}^3$, the attitude matrices of the two quadrotors, 
$(R_1, R_2)\in (SO(3))^2$ and the directions of the links which connect the center of mass of each quadrotor respectively with the mass point, 
$(q_1,q_2)\in (S^2)^2$. The configuration manifold of the system is $Q=\mathbb{R}^3\times (SO(3))^2 \times (S^2)^2$.
\begin{figure}[htbp]
    \centering
    \includegraphics[width=.9\textwidth]{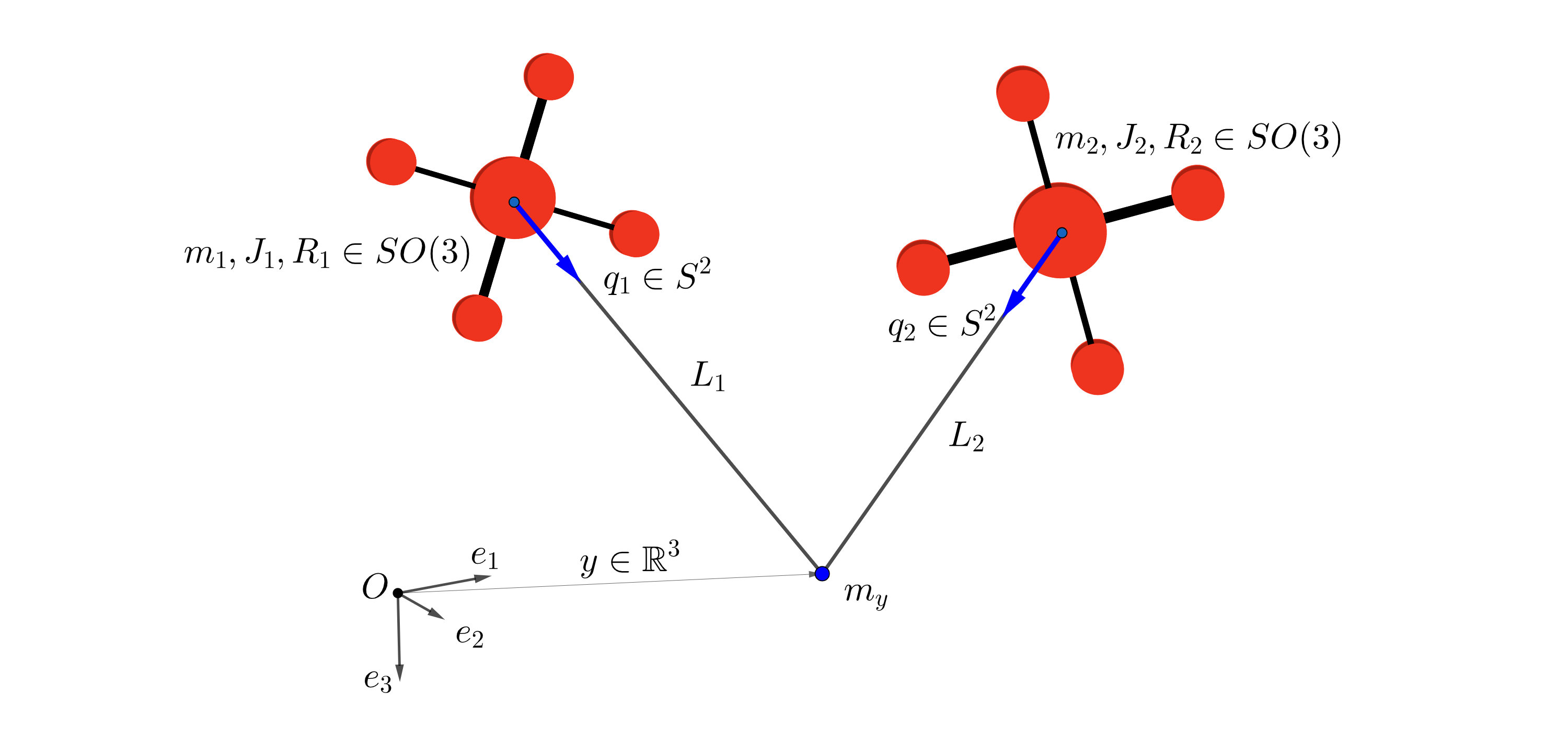}
    \caption{Two quadrotors connected to the mass point $m_y$ via massless links of lengths $L_i$.}
    \label{fig:scheme}
\end{figure}
In order to present the equations of motion of the system we start by identifying $TSO(3)\simeq SO(3)\times \mathfrak{so}(3)$ via left-trivialization. This choice allows us to write the kinematic equations of the system as 
\begin{equation}\label{eq:kinematic}
\dot{R}_i = R_i\hat{\Omega}_i,\quad \dot{q}_i = \hat{\omega}_iq_i\quad \quad i=1,2,
\end{equation} where $\Omega_1,\Omega_2\in\mathbb{R}^3$ represent the angular velocities of each quadrotor, respectively, and $\omega_1,\omega_2$ express the time derivatives of the orientations $q_1,q_2\in S^2$, respectively, in terms of angular velocities, expressed with respect to the body-fixed frames. From these equations we define the trivialized Lagrangian \[
L(y,\dot{y},R_1,\Omega_1,R_2,\Omega_2,q_1,\omega_1,q_2,\omega_2): \mathbb{R}^6\times \left(SO(3)\times \mathfrak{so}(3)\right)^2\times (TS^2)^2\rightarrow \mathbb{R},
\] as the difference of the total kinetic energy of the system and the total potential (gravitational) energy, $L=T-U$, with:

\[T = \frac{1}{2}m_y\|\dot{y}\|^2 +\frac{1}{2}\sum_{i=1}^2 (m_i\|\dot{y} -L_i\hat{\omega}_iq_i \|^2 + \Omega_i^TJ_i\Omega_i) ,\]

and $$U= -m_yge_3^Ty - \sum_{i=1}^2 m_ige_3^T(y-L_iq_i),$$
where $J_1,J_2\in\mathbb{R}^{3\times 3}$ are the inertia matrices of the two quadrotors and $m_1,m_2\in\mathbb{R}^+$  are their respective total masses. In this system each of the two quadrotors generates a thrust force, which we denote with $u_i = -T_iR_ie_3\in\mathbb{R}^3$, where $T_i$ is the magnitude, while $e_3$ is the direction of this vector in the $i-$th body-fixed frame, $i=1,2$. The presence of these forces make it a non conservative system. Moreover, the rotors of the two quadrotors generate a moment vector, and we denote with $M_1, M_2\in\mathbb{R}^3$ the cumulative moment vector of each of the two quadrotors. To derive the Euler--Lagrange equations, a possible approach is through Lagrange--d'Alambert's principle, as presented in \cite{lee18gfo}. We write them in matrix form as
\begin{equation}  \label{systquad}
A(z)\dot{z} = h(z)
\end{equation}
where
\[
z = [y,v,\Omega_1,\Omega_2,\omega_1,\omega_2]^T\in\mathbb{R}^{18},\] 
\[
A(z) = \begin{bmatrix} I_3 & 0_3 & 0_3 & 0_3 & 0_3 & 0_3 \\ 0_3 & M_q  & 0_3 & 0_3  & 0_3 & 0_3   \\ 0_3 & 0_3 & J_1 & 0_3 & 0_3 & 0_3 \\ 0_3 & 0_3 & 0_3 & J_2 &  0_3 &  0_3 \\ 0_3 & -\frac{1}{L_1}\hat{q}_1 & 0_3 & 0_3 & I_3 & 0_3 \\ 0_3 & -\frac{1}{L_2}\hat{q}_2 & 0_3 & 0_3 & 0_3 & I_3\end{bmatrix},
\]
\[
h(z) = \begin{bmatrix}h_1(z) \\ h_2(z) \\ h_3(z) \\ h_4(z) \\  h_5(z) \\ h_6(z)\end{bmatrix} =\begin{bmatrix} v \\ -\sum_{i=1}^{2} m_{i}L_{i}\|\omega_{i}  \|^{2} q_{i} + M_q g e_{3}+\sum_{i=1}^{2} u_i^{\parallel} \\ -\Omega_1\times J_1\Omega_1 + M_1 \\ -\Omega_2\times J_2\Omega_2 + M_2 \\ -\frac{1}{L_1} g \hat{q}_{1} e_{3} -\frac{1}{m_1L_1}q_{1} \times u_1^{\perp}\\ -\frac{1}{L_2} g \hat{q}_{2} e_{3} -\frac{1}{m_2L_2}q_{2} \times u_2^{\perp}\end{bmatrix},
\]
where $M_q = m_yI_3 + \sum_{i=1}^2m_iq_iq_i^T,$ and  $u_i^{\parallel},u_i^{\perp}$ are respectively the orthogonal projection of $u_i$ along $q_i$ and to the plane $T_{q_i}S^2$, $i=1,2$, i.e. $u_i^{\parallel}=q_{i} q_{i}^{T}u_i$, $u_i^{\perp}=(I-q_{i} q_{i}^{T})u_i$. 
These equations, coupled with the kinematic equations in (\ref{eq:kinematic}), describe the dynamics of a point 
$$
P = \left[y ,\;\; v,\;\; R_1 ,\;\; \Omega_1 ,\;\; R_2 ,\;\; \Omega_2 ,\;\; q_1 ,\;\; \omega_1  ,\;\; q_2 ,\;\; \omega_2 \right] \in M = TQ.
$$
Since the matrix $A(z)$ is invertible, we pass to the following set of equations
\begin{equation} \label{eq:dynamics} 
\dot{z} = A^{-1}(z)h(z):=\Tilde{h}(z) :=\bar{h}(P) = [\bar{h}_1(P),...,\bar{h}_7(P)]^T.
\end{equation}
\subsection{Analysis via transitive group actions}
We identify the phase space M with $M\simeq T\mathbb{R}^3\times (TSO(3))^2 \times (TS^2)^2$. The group we consider is
$$\bar{G} = \mathbb{R}^6 \times (TSO(3))^2 \times (SE(3))^2,$$ where the groups are combined with a direct-product structure and $\mathbb{R}^6$ is the additive group.
For a group element $$g=((a_1,a_2),((B_1,b_1),(B_2,b_2)),((C_1,c_1),(C_2,c_2)))\in \bar{G}$$ and a point $P \in M$ in the manifold,
we consider the following left action

\[
\begin{split}
\psi_g(P) = [y+a_1, \;\;v+a_2,\;\; &B_1R_1,\;\;  \Omega_1 + b_1,\;\; B_2R_2,\;\; \Omega_2 + b_2,\;\;\\ &C_1q_1,\;\;C_1\omega_1 + c_1\times C_1q_1,\;\; C_2q_2,\;\;C_2\omega_2 + c_2\times C_2q_2
].
\end{split}
\]
The well-definiteness and transitivity of this action come from standard arguments, see for example \cite{olver2000applications}. 
The infinitesimal generator associated to $$\xi = \left[\xi_1 ,\;\; \xi_2,\;\; \eta_1 ,\;\; \eta_2 ,\;\; \eta_3 ,\;\; \eta_4 ,\;\; \mu_1 ,\;\; \mu_2 ,\;\; \mu_3 ,\;\; \mu_4 \right]\in \mathfrak{\bar{g}},$$ where $ \mathfrak{\bar{g}}=T_e\bar{G}$, writes
\[\begin{split}
\infgen(\xi)\vert_P = [\xi_1,\;\; \xi_2, \;\; \hat{\eta}_1R_1,\;\; \eta_2,\;\; &\hat{\eta}_3R_2,\;\;  \eta_4,\;\;\\ & \hat{\mu}_1q_1,\;\; \hat{\mu}_1\omega_1 + \hat{\mu}_2q_1, \;\; \hat{\mu}_3q_2,\;\; \hat{\mu}_3\omega_2 + \hat{\mu}_4q_2 ].
\end{split}
\]
We can now focus on the construction of the function $f:M\rightarrow \bar{\mathfrak{g}}$ such that $\infgen(f(P))\vert_P=F\vert_P$, where

\[\begin{split}
F\vert_P = [\bar{h}_1(P), \;\; \bar{h}_2(P), \;\; R_1&\hat{\Omega}_1,\;\; \bar{h}_3(P), \;\;  R_2\hat{\Omega}_2,\;\;\\  &\bar{h}_4(P), \;\; \hat{\omega}_1q_1, \;\; \bar{h}_5(P),\;\; \hat{\omega}_2q_2, \;\; \bar{h}_6(P)
]\in T_{P}M
\end{split}\]
is the vector field obtained combining the equations (\ref{eq:kinematic}) and (\ref{eq:dynamics}). 
We have
\[
\begin{split}
f(P) = [\bar{h}_1(P),\;\; \bar{h}_2(P),\;\; R_1\Omega_1,\;\;&\bar{h}_3(P),\;\; R_2\Omega_2,\;\;\bar{h}_4(P),\\ \;\;&\omega_1,\;\; q_1\times \bar{h}_5(P),\;\;\omega_2,\;\; q_2\times \bar{h}_6(P)]\in\bar{\mathfrak{g}}.
\end{split}
\]
We have obtained the  local representation of the vector field $F\in\mathfrak{X}(M)$ in terms of the infinitesimal generator of the transitive group action $\psi$, hence we can solve for one time step $\Delta t$ the IVP
\[
\begin{cases}
\dot{\sigma}(t) = \dexp_{\sigma(t)}^{-1}\Big(f\big(\psi(\exp(\sigma(t)),P(t))\big)\Big) \\
\sigma(0) = 0\in\bar{\mathfrak{g}}
\end{cases}
\]
and then update the solution $P(t+\Delta t) = \psi(\exp(\sigma(\Delta t)),P(t))$. 

The above construction is completely independent of the control functions $\{u_i^{\parallel},u_i^{\perp},M_i\}_{i=1,2}$ and hence it is compatible with any choice of these parameters.

\subsection{Numerical experiments}

We tested Lie group numerical integrators for a load transportation problem presented in \cite{Lee2013GeometricCO}. The control inputs $\{u_i^{\parallel},u_i^{\perp},M_i\}_{i=1,2}$ are constructed such that the
point mass asymptotically follows a given desired trajectory $y_d \in \mathbb{R}^3$, given by a smooth function of time,  and the
quadrotors maintain a prescribed formation relative to
the point mass. In particular, the parallel components $u_i^{\parallel}$
are designed such that the payload follows the desired trajectory $y_d$ (load transportation problem), while the normal components $u_i^{\perp}$ are designed such that  $q_i$ converge to desired directions $q_{id}$ (tracking problem in $S_2$). Finally, ${M_i}$ are designed to control the attitude of the quadrotors.

In this experiment we focus on a simplified dynamics model, i.e. we neglect the construction of the controllers $M_i$ for the attitude dynamics of the quadrotors. However, the full dynamics model can also be easily integrated, once the expressions for the attitude controllers are available.

\begin{figure}[ht]
    \centering
    \includegraphics[trim= 200 0 300 0,width=.9\textwidth]{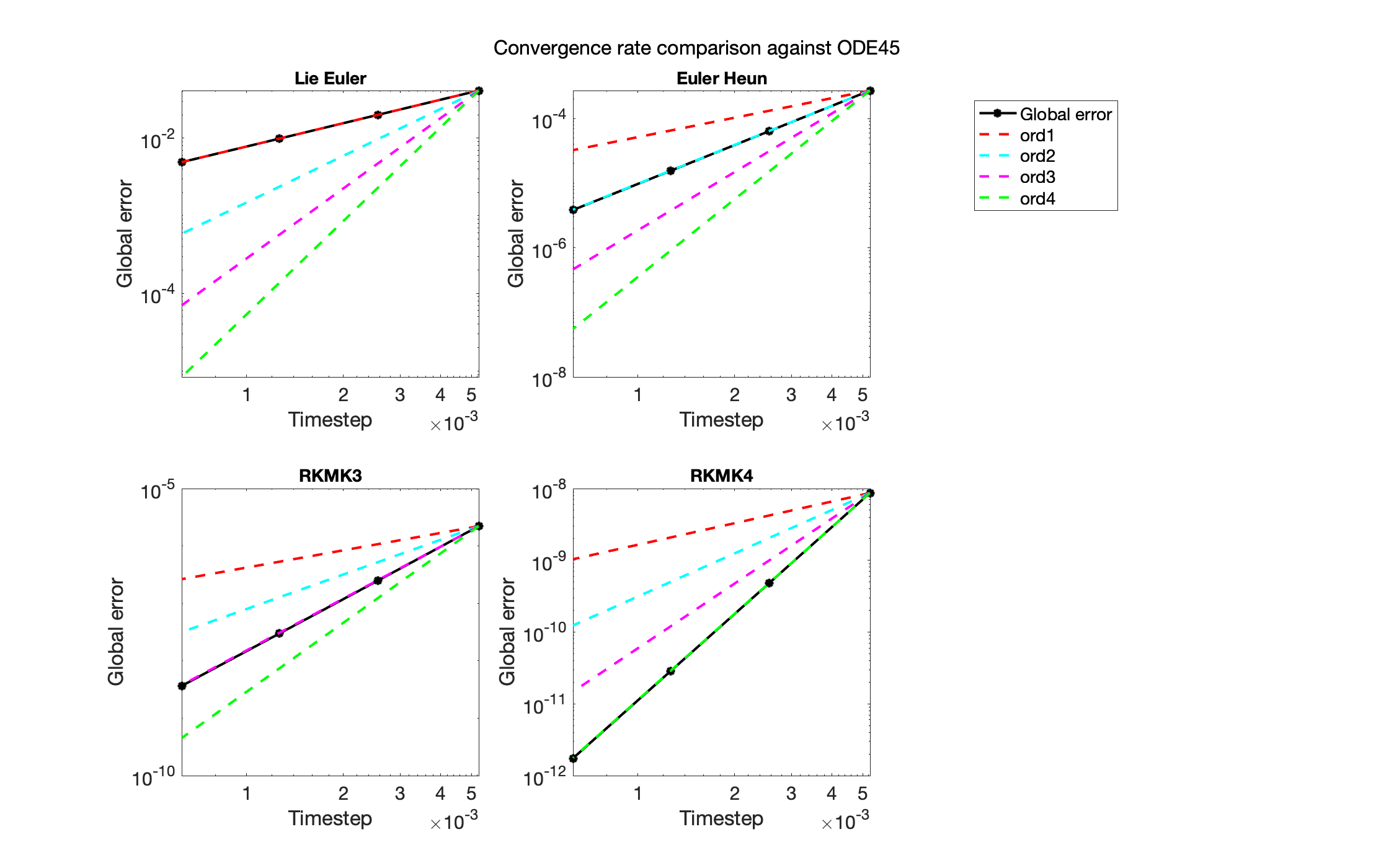}
    \caption{Convergence rate of the numerical schemes compared with ODE45}
    \label{fig:convRatesQuad}
\end{figure}

In Figure \ref{fig:convRatesQuad} we show the convergence rate of four different RKMK methods \RE{compared with the reference solution obtained with ODE45 in MATLAB}.

\begin{figure}[htbp]
  \centering
  \begin{minipage}[b]{0.495\textwidth}
    \includegraphics[trim=190 110 190 110,clip,width=\textwidth]{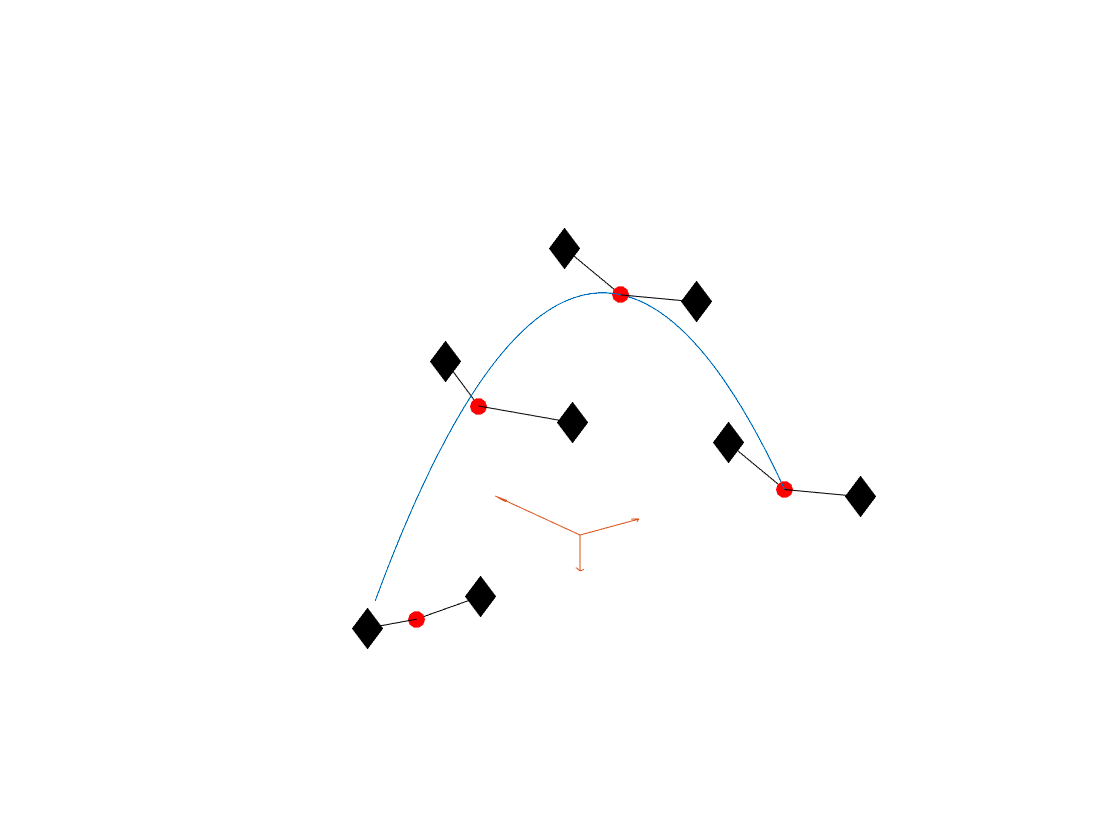}
    \caption{Snapshots at $0 \leq t \leq 5$. 
    \newline
    \newline}   
    \label{fig1quad}
  \end{minipage}
  \begin{minipage}[b]{0.495\textwidth}
    \includegraphics[trim=30 0 30  0,clip,width=\textwidth]{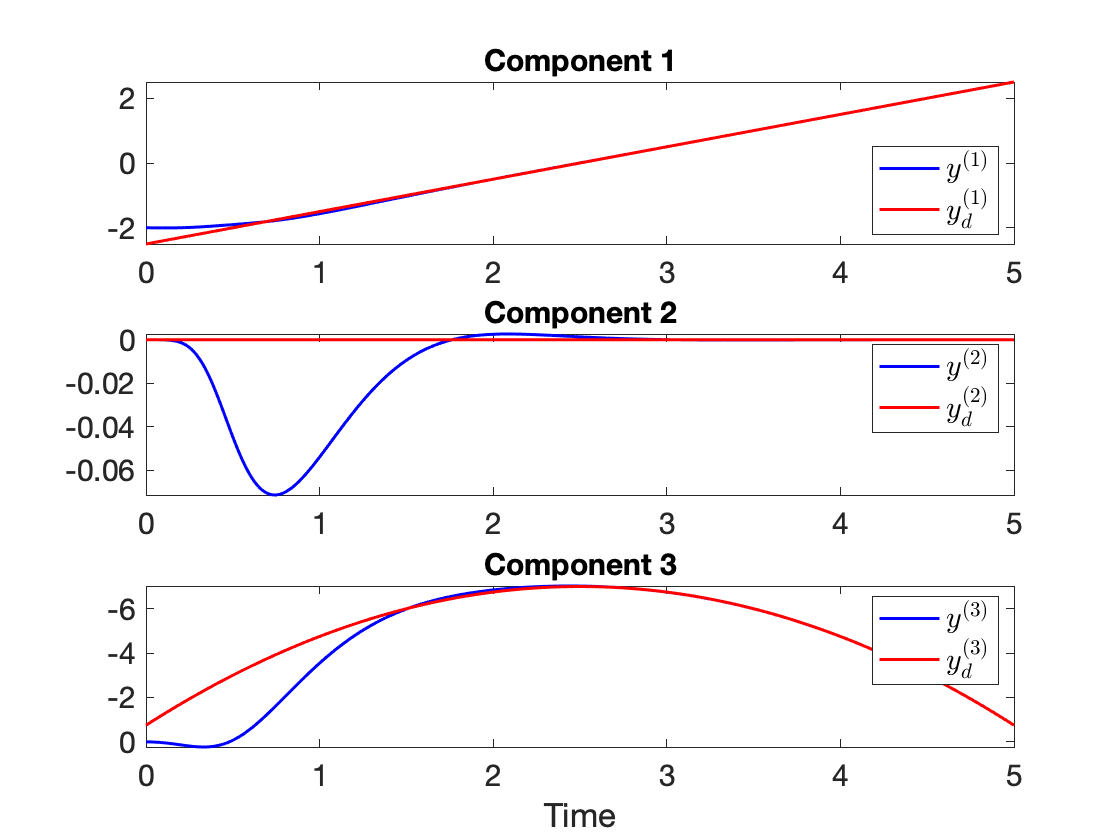}
    \caption{Components of the load position (in blue) and the desired trajectory (in red) as a function time.}
     \label{fig2quad}
  \end{minipage}
\end{figure}

\begin{figure}[htbp]
  \centering
  \begin{minipage}[b]{0.495\textwidth}
    \includegraphics[trim= 25 0 30 20,clip,width=\textwidth]{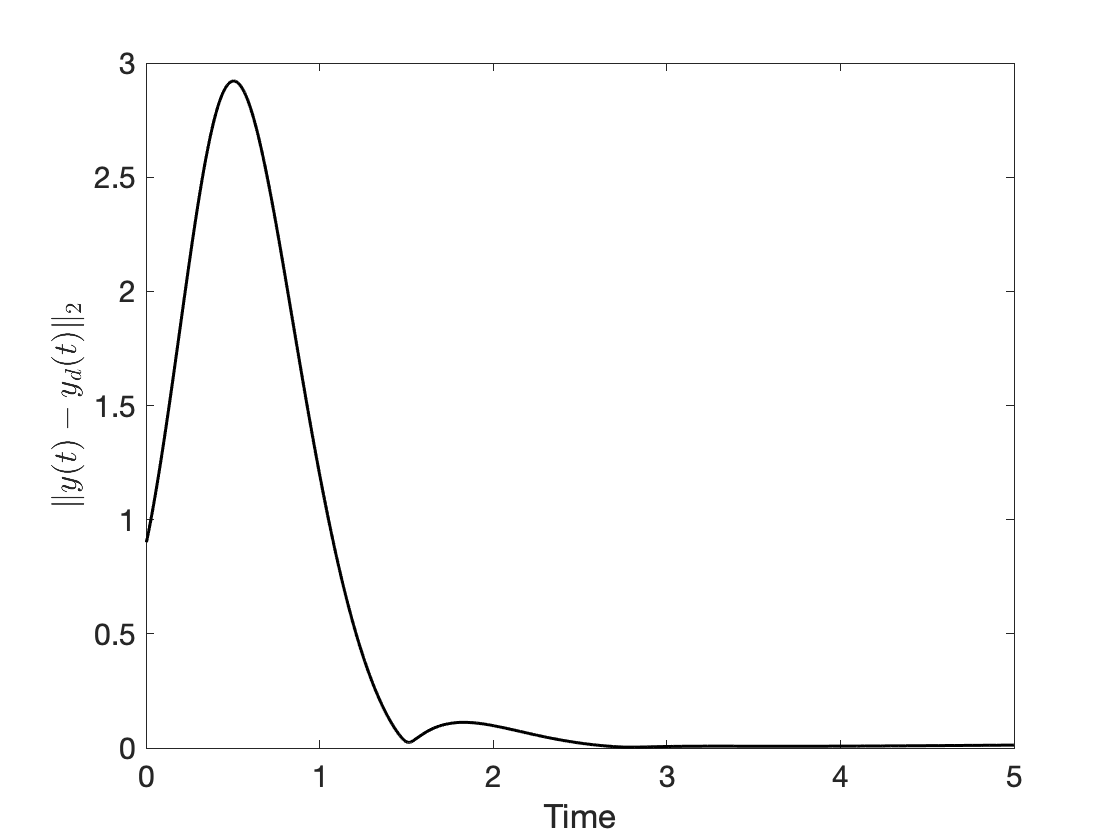}
    \caption{Deviation of the load position from the target trajectory. }
    \label{fig3quad}
  \end{minipage}
  \begin{minipage}[b]{0.495\textwidth}
    \includegraphics[trim= 20 0 30 0,clip,width=\textwidth]{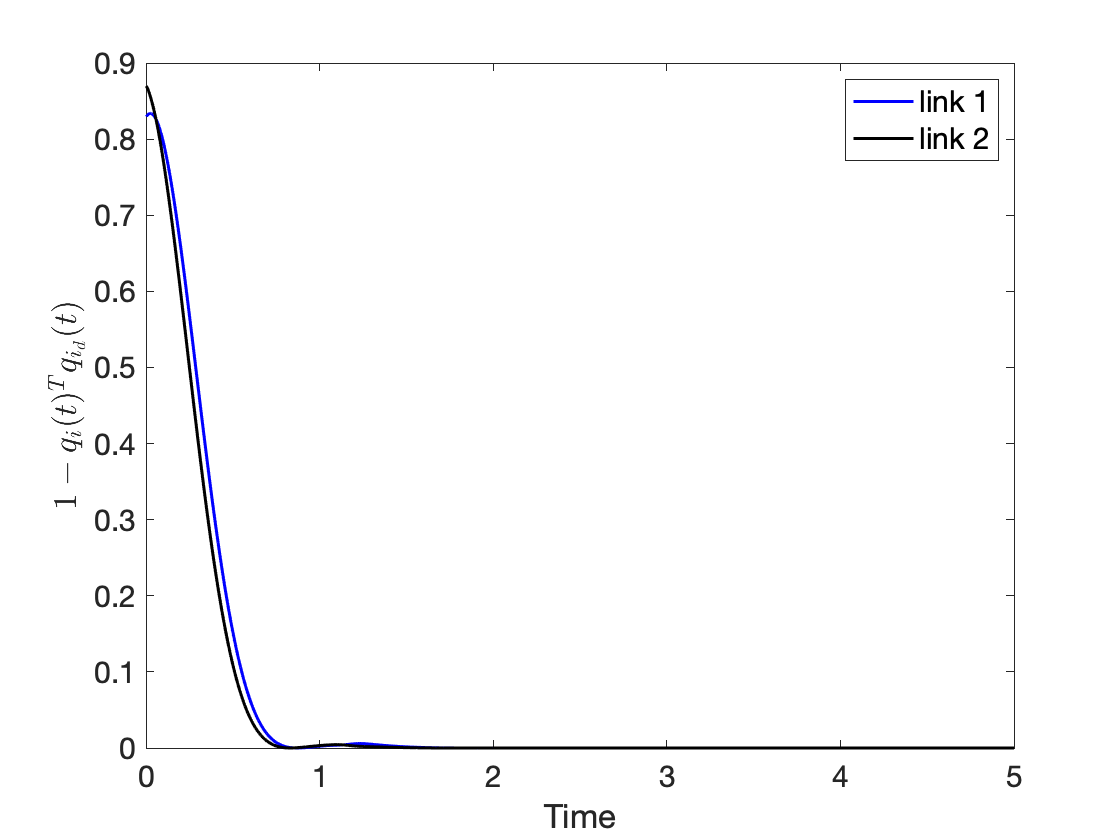}
    \caption{Direction error of the links.
    \newline}
     \label{fig4quad}
  \end{minipage}
\end{figure}

\begin{figure}[htbp]
    \centering
    \includegraphics[trim = 50 0 50 0, clip, width=1\textwidth]{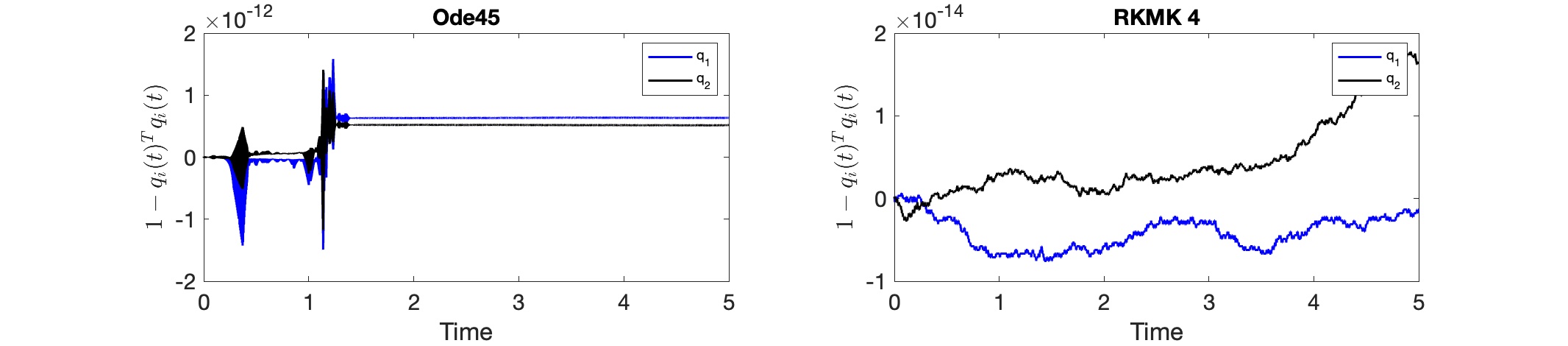}
    \caption{Preservation of the norms of $q_1, q_2 \in S^2$.}
    \label{fig5quad}
\end{figure}

In Figures \ref{fig1quad}-\ref{fig5quad} we show results in the tracking of a parabolic trajectory, obtained by integrating the system \eqref{systquad} with a RKMK method of order 4. 

\section{Summary and outlook}
\RE{In this paper we have considered Lie group integrators with a particular focus on problems from mechanics. In mathematical terms this means that the Lie groups and manifolds of particular interest are $SO(n),\ n=2,3$, $SE(n),\ n=2,3$ as well as the manifolds $S^2$ and $TS^2$. 
The abstract formulations by e.g. Crouch and Grossman \cite{crouch93nio}, Munthe-Kaas \cite{munthe-kaas99hor} and Celledoni et al. \cite{celledoni03cfl} have often been demonstrated on small toy problems in the literature, such as the free rigid body or the heavy top systems. But in papers like \cite{bruls10otu}, hybrid versions of Lie group integrators have been applied to more complex beam and multi-body problems. The present paper is attempting to move in the direction of more relevant examples without causing the numerical solution to depend on how the manifold is embedded in an ambient space, or the choice of local coordinates.}

\RE{It will be the subject of future work to explore more examples and to aim for a more systematic approach to applying Lie group integrators to mechanical problems. 
In particular, it is of interest to the authors to consider models of beams, that could be seen as a generalisation of the $N$-fold pendulum discussed here.}

\clearpage 

\end{document}